\documentclass{amsart}

\usepackage{a4wide, amssymb, mathrsfs, stmaryrd}
\usepackage{fancyheadings}

\usepackage[
open,
openlevel=2,
atend
]{bookmark}[2009/08/13]

\usepackage{amsmath}

\usepackage[T1]{fontenc}

\usepackage[arrow,matrix,curve]{xy}

\usepackage{color}

\def\xto#1{\xrightarrow{#1}}

\newcommand{\nc}{\newcommand}

\newcommand{\myMO}[1]{{\fontshape{rm}{\textbf{#1}}}}

\DeclareMathOperator{\Grmod}{\myMO{Grmod}\hspace{+0.25ex}-\hspace{-0.25ex}}
\DeclareMathOperator{\Mod}{\myMO{Mod}\hspace{+0.25ex}-\hspace{-0.25ex}}
\DeclareMathOperator{\proj}{\myMO{proj.dim}}
\DeclareMathOperator{\gl}{\myMO{gl.dim}}

\newcommand*\cocolon{%
        \nobreak
        \mskip6mu plus1mu
        \mathpunct{}%
        \nonscript
        \mkern-\thinmuskip
        {:}%
        \mskip2mu
        \relax
}

\nc{\Ch}{\operatorname{Ch}}
\nc{\sA}{{\mathscr A}}
\nc{\wt}{\widetilde} \nc{\bl}{\bullet} \nc{\al}{\alpha}
\nc{\sg}{\sigma} \nc{\vf}{\varphi} \nc{\om}{\omega}
\nc{\ve}{\varepsilon} \nc{\ol}{\overline} \nc{\lb}{\lambda}
\nc{\Lb}{\Lambda} \nc{\Gm}{\Gamma} \nc{\cP}{{\mathscr P}}
\nc{\sB}{{\mathscr B}}
\nc{\ul}{\underline} \nc{\os}{\overset} \nc{\us}{\underset}
\nc{\pa}{\partial} \nc{\wh}{\widehat} \nc{\sbs}{\subset} \nc{\br}{\breve}
\nc{\lra}{\longrightarrow} \nc{\all}{\allowdisplaybreaks}
\nc{\Ker}{\operatorname{Ker}} \nc{\Img}{\operatorname{Im}}
\nc{\Kan}{\operatorname{Kan}} \nc{\Hom}{\operatorname{Hom}}
\nc{\Imm}{\operatorname{Im}}   \nc{\Ho}{\operatorname{Ho}}
\nc{\Ext}{\operatorname{Ext}}    \nc{\Cone}{\operatorname{Cone}}
\nc{\pr}{\operatorname{pr}} \nc{\cls}{\operatorname{cls}}
\nc{\cof}{\operatorname{cof}}
\nc{\sSet}{\operatorname{\textbf{sSet}}}
\nc{\Map}{\operatorname{Map}}
\nc{\incl}{\operatorname{incl}}
\nc{\Hocolim}{\operatorname{Hocolim}}
\nc{\colim}{\operatorname{colim}}
\nc{\Endd}{\operatorname{End}}
\nc{\const}{\operatorname{const}}
\nc{\inn}{\operatorname{in}}
\nc{\Ev}{\operatorname{Ev}}
\nc{\rr}{\operatorname{\textbf{r}}}
\nc{\cop}{\operatorname{\textbf{l}}}

\numberwithin{equation}{subsection}
\newtheorem{theo}[equation]{Theorem}
\newtheorem{lem}[equation]{Lemma}
\newtheorem{prop}[equation]{Proposition}
\newtheorem{coro}[equation]{Corollary}
\theoremstyle{definition}
\newtheorem{defi}[equation]{Definition}
\newtheorem{conv}[equation]{Convention}
\newtheorem{remk}[equation]{Remark}
\newtheorem{exmp}[equation]{Example}
\swapnumbers
\newtheorem{numbered paragraph}[equation]{}

\newcommand{\mybox}{\ensuremath \Box}

\newenvironment{prf}{\noindent\textbf{Proof.}}{\mybox}

\begin{document}

\def\S{\textbf{S}}
\def\Z{{\mathbb Z}}
\def\L{{\mathcal L}}
\def\B{{\mathcal B}}
\def\Q{{\mathcal Q}}
\def\M{{\mathcal M}}
\def\D{{\mathcal D}}
\def\R{{\mathscr R}}
\def\E{{\mathcal E}}
\def\K{{\mathcal K}}
\def\W{{\mathcal W}}
\def\N{{\mathcal N}}
\def\T{{\mathcal T}}
\def\A{{\mathcal A}}
\def\C{{\mathcal C}}
\def\P{{\mathcal P}}
\def\V{{\mathcal V}}
\def\G{{\mathcal G}}
\let\xto\xrightarrow

\title[]
{On the algebraic classification of module spectra}

\author{Irakli Patchkoria}

\thanks{This research was supported by the Deutsche Forschungsgemeinschaft Graduirtenkolleg 1150 ``Homotopy and Cohomology''}

\begin{abstract}
Using methods developed by Franke in \cite{F96}, we obtain algebraic classification results for modules over certain symmetric ring spectra ($S$-algebras).
In particular, for any symmetric ring spectrum $R$ whose graded homotopy ring $\pi_*R$ has graded global homological dimension $2$ and
is concentrated in degrees divisible by some natural number $N \geq 4$, we prove that the homotopy category of $R$-modules is equivalent
to the derived category of the homotopy ring $\pi_*R$. This improves the Bousfield-Wolbert algebraic classification of isomorphism classes of objects of
the homotopy category of $R$-modules. The main examples of ring spectra to which our result applies are
the $p$-local real connective $K$-theory spectrum $ko_{(p)}$, the Johnson-Wilson spectrum $E(2)$, and the truncated Brown-Peterson spectrum $BP\langle 1 \rangle$, for an odd prime $p$.
\end{abstract}


\maketitle

\setcounter{section}{0}

\section{Introduction}

\setcounter{subsection}{1}

In \cite{B85} Bousfield gave an algebraic classification of isomorphism classes of objects in the $E(1)$-local (or, equivalently, $K_{(p)}$-local) stable homotopy category $L_1 \S$ for an odd prime $p$. For this he used a certain $k$-invariant coming from a $d_2$-differential of the $E(1)$-homology Adams spectral sequence. Bousfield also tried to describe morphisms in $L_1 \S$ but only managed to characterize them up to Adams filtration.

In \cite{F96} Franke, using the language of triangulated derivators, proved an abstract uniqueness theorem for triangulated categories possessing an Adams spectral sequence. Applying this result to $L_1 \S$ (for an odd prime), he obtained an equivalence of categories between $L_1 \S$ and the derived category of certain twisted chain complexes of $E(1)_*E(1)$-comodules. (For a streamlined exposition of this result see \cite{R08}.) This gives a complete algebraic description of the category $L_1 \S$ and thus improves the aforementioned result of Bousfield. Moreover, as is shown in \cite{R08}, this equivalence is exotic in the sense that it does not come from a zig-zag of Quillen equivalences.

A similar classification result to that of Bousfield for module spectra was obtained by Wolbert in \cite{W98}. More precisely, let $R$ be a symmetric ring spectrum (or an $S$-algebra) and suppose that the (right) graded global homological dimension of the homotopy ring $\pi_*R$ is less or equal than two. Suppose further that $\pi_*R$ is concentrated in dimensions divisible by some natural number $N \geq 2$. Under these assumptions, Wolbert gives an algebraic characterization of isomorphism classes of objects in the derived category $\D(R)$ \cite[Theorem 6]{W98}. The main example of a ring spectrum to which this can be applied is the $p$-local real connective $K$-theory spectrum $ko_{(p)}$ for $p$ an odd prime.

Note that Wolbert, like Bousfield, did not give a complete algebraic characterization of morphisms in $\D(R)$. Consequently, a natural question arises whether one can apply Franke's methods to improve Wolbert's classification. One of the aims of the present paper is to examine this question.

Let $R$ be a symmetric ring spectrum. We say that an object $X$ of $\D(R)$ has dimension $k$ if the projective dimension of the homotopy $\pi_*X$ regarded as a graded $\pi_*R$-module is equal to $k$. Similarly, we say that an object $M$ in the derived category $\D(\pi_*R)$ (i.e., a differential graded $\pi_*R$-module) has dimension $k$ if the projective dimension of the homology $H_*M$ regarded as a graded $\pi_*R$-module is equal to $k$. The main aim of this paper is to give a self-contained homotopy theoretic proof of the following theorems:

\begin{theo}\label{Rcoro} Let $R$ be a symmetric ring spectrum. Suppose that the graded homotopy ring $\pi_*R$ is concentrated in dimensions divisible by a natural number $N$ and assume that the (right) graded global homological dimension of $\pi_*R$ is less than $N$. Then one can construct a functor
$$\R \colon \D(\pi_*R) \longrightarrow \D(R)$$
such that the following hold:

{\rm (i)} The diagram
$$\xymatrix{\D(\pi_*R) \ar[rr]^{\R} \ar[dr]_{H_*} & & \D(R) \ar[dl]^{\pi_*} \\ & \Grmod{\pi_*R}, & }$$
where $\Grmod{\pi_*R}$ denotes the category of (right) graded $\pi_*R$-modules, commutes up to a natural isomorphism.

{\rm (ii)} The functor $\R \colon \D(\pi_*R) \longrightarrow \D(R)$ restricts to an equivalence of the full subcategories of at most one-dimensional objects.

{\rm (iii)} Let $X \in \D(R)$ and suppose that $X$ has dimension at most two. Then there exists $M$  in $\D(\pi_*R)$ such that $\R(M)$ is isomorphic to $X$
in $\D(R)$. \end{theo}

\begin{theo}\label{equival1} Let $R$ be a symmetric ring spectrum. Suppose that the graded homotopy ring $\pi_*R$ is concentrated in dimensions divisible by a natural number $N$ and assume that the (right) graded global homological dimension of $\pi_*R$ is less than $N-1$. Then the functor $$\R \colon \D(\pi_*R) \longrightarrow \D(R)$$ restricts to an equivalence of the full subcategories of at most two-dimensional objects. \end{theo}

In particular, we have

\begin{coro}\label{equival2} Let $R$ be a symmetric ring spectrum. Suppose that the graded homotopy ring $\pi_*R$ is concentrated in dimensions divisible by a natural number $N \geq 4$ and assume that the (right) graded global homological dimension of $\pi_*R$ is equal to two. Then the functor $$\R \colon \D(\pi_*R) \longrightarrow \D(R)$$ is an equivalence of categories.
\end{coro}

This corollary improves Wolbert's classification (\cite[Theorem 6]{W98}) for $N \geq 4$.

Note that similar results hold for differential graded rings and our proofs work in this algebraic setting as well (see Remark \ref{algebra}).

The functor $\R \colon \D(\pi_*R) \longrightarrow \D(R)$ is in fact Franke's functor for the special case of $\D(R)$. Franke constructs his functor for general triangulated derivators possessing an Adams spectral sequence, and claims that it is an equivalence of categories \cite[2.2, Proposition 2]{F96}. If this is so, then under the hypotheses of \ref{Rcoro} the functor $\R \colon \D(\pi_*R) \longrightarrow \D(R)$ becomes an equivalence that trivially implies \ref{Rcoro} and \ref{equival1}. However, as far as the present author can see, the proof of \cite[2.2, Proposition 2]{F96} seems to contain a gap (see Remark \ref{shecdoma}) which we were unable to fill in the general setting. However, we managed to eliminate this gap in low dimensional cases and thus obtained the statements of \ref{Rcoro} and \ref{equival1}.

Note that for the construction of $\R$ we adopt the methods of \cite{R08}. In particular, we use the language of model categories rather than that of derivators which is used by Franke. This technically simplifies our exposition.

For $n <  2(p-1)$, the Johnson-Wilson spectra $E(n)$ are non-trivial examples of ring spectra to which Theorem \ref{Rcoro} can be applied. Further important examples are the truncated Brown-Peterson spectra $BP\langle n \rangle$ for $ n+1 < 2(p-1)$. (For all primes and any $n$, $E(n)$ as well as $BP\langle n \rangle$ admit $A_\infty$-structures (see e.g. \cite{L03}) and thus possess models which are symmetric ring spectra.) Note that for all these ring spectra the functor $\R$ is not a derived functor of a Quillen functor (see \ref{Rnonderived}).

Examples of ring spectra which fulfill the hypotheses of \ref{equival2} are the $p$-local real connective $K$-theory spectrum $ko_{(p)}$, the Johnson-Wilson spectrum $E(2)$, and the truncated Brown-Peterson spectrum $BP\langle 1 \rangle$, for an odd prime $p$. We show in \ref{Rnonderived} that all these examples are exotic in the sense that the equivalence $\R$ does not come from a zig-zag of Quillen equivalences.

The paper is organized as follows. In Section \ref{prele} some basic facts about triangulated categories and model categories are recalled. Next we briefly review stable model categories and give examples which are important for our purposes. Furthermore we recall model structures on diagram categories and long exact sequences of a homotopy pullback and homotopy coequalizer.

In Section \ref{yv. mtavari} the functor $\R$ is constructed and \ref{Rcoro} (i) is proved. In fact, we define $\R$ in a more general setting of stable model categories. The most complicated step in the construction of $\R$ is discussed in Subsection \ref{Q} which concerns Franke's crowned diagrams (see \ref{exmpV}). At the end of Section \ref{yv. mtavari} we verify that the functor $\R$ commutes with suspensions.

In Section \ref{appl} we recall the notion of a derived mapping space and prove a proposition which gives a sufficient condition for  $\R$ not being a derived functor. We also discuss some non-trivial examples to which \ref{Rcoro} applies.

Section \ref{advgamoy} is devoted to the proof of \ref{Rcoro} (ii). As a corollary we get that under the hypothesis of \ref{Rcoro} the functor $\R \colon \D(\pi_*R) \longrightarrow \D(R)$ is an equivalence of categories if the graded global homological dimension of $\pi_*R$ is equal to one. The complex $K$-theory spectrum $KU$ and the connective Morava $K$-theory spectra $k(n)$, $n \geq 1$ serve as our main examples here.

Finally, in Section \ref{intgamoy} we prove \ref{Rcoro} (iii) and \ref{equival1}. We also discuss non-trivial examples to which \ref{equival2} applies. Note that the key results in this section are \ref{techlemma} and \ref{dreieck1} which fill the aforementioned gap in the proof of \cite[2.2, Proposition 2]{F96} for the two dimensional case.

\section*{Acknowledgments}

This paper is based on my master's thesis at the University of Bonn. I would like to thank my advisor Stefan Schwede for suggesting this topic and for his help and encouragement throughout the way. Many thanks also go to Katja Hutschenreuter and Lennart Meier for their valuable comments and suggestions on earlier drafts of the paper.

\section{Preliminaries}\label{prele}

\subsection{Notation}

\subsubsection{Simplicial objects}\label{sset}

The category of simplicial sets is denoted by \textbf{sSet}. For $n \geq 0$, $\Delta[n]$ stands for the standard simplicial $n$-simplex. Next, if $K$ is a simplicial set and $x \in K_n$, then we denote again by $x$ the unique simplicial map $\Delta[n] \longrightarrow K$ which sends the generator of $\Delta[n]$ to $x$.

The letter $I$ stands for $\Delta[1]$ and $I_+$ for the union of $\Delta[1]$ with a disjoint base point. $S^0$ denotes the pointed simplicial $0$-sphere, i.e., the union of the standard $0$-simplex with a disjoint base point. Further,
$$\inn_0, \inn_1 \colon S^0 \longrightarrow I_+$$
are our notations for obvious maps sending the non-base point of $S^0$ to $(0)$ and $(1)$, respectively. $S^1$ denotes the simplicial circle $I/(0 \sim 1)$, i.e., $\Delta[1]/\partial \Delta[1]$.

\subsubsection{DG modules and graded modules} \label{grdims} For any differential graded ring $A$, we denote by $\Mod A$ the category of differential graded right $A$-modules. In particular, for any graded ring $B$, we have the category $\Mod B$, where $B$ is regarded as a differential graded ring with zero differential. Further, for any graded ring $B$, $\Grmod B$ stands for the category of graded right modules over $B$. For any $X \in \Grmod B$, $\proj X$ denotes the graded projective dimension of $X$, and $\gl B$ the (right) graded global homological dimension of $B$.

\subsubsection{Morphism objects} If $\T$ is a triangulated category and $X, Y \in \T$, then $\Hom_{\T}(X,Y)_*$ is our notation for $\Hom_{\T}(\Sigma^*X,Y)$,
the graded Hom. Next, for a stable model category $\M$ and objects $X, Y \in \M$, the abelian group of morphisms from $X$ to $Y$ in the homotopy
category is denoted by $[X,Y]$. (The notions of a triangulated category and a model category will be recalled in the following subsections.)

\vspace{0.75em}

\subsection{Triangulated categories}~

\vspace{0.75em}

\noindent Let $\T$ be a category equipped with an endofunctor $\Sigma \colon \T \longrightarrow \T$. A triangle in $\T$ with respect to $\Sigma$ is a  sequence
$$\xymatrix{ X \ar[r]^{f} & Y \ar[r]^{g}  & Z \ar[r]^{h} & \Sigma X.}$$
A morphism $(\alpha, \beta, \gamma) \colon (f,g,h) \longrightarrow (f',g',h')$ of triangles is a commutative diagram
$$\xymatrix{ X \ar[r]^f \ar[d]_{\alpha} & Y \ar[r]^g \ar[d]^{\beta} & Z \ar[r]^h \ar[d]^{\gamma} & \Sigma X \ar[d]^{\Sigma \alpha} \\ X' \ar[r]^{f'} & Y' \ar[r]^{g'} & Z \ar[r]^{h'}  & \Sigma X'.}$$
A morphism $(\alpha, \beta, \gamma)$  is an isomorphism of triangles if all three components are isomorphisms.

\begin{defi}\label{def0} A \emph{triangulated category} is an additive category $\T$ together with a self-equivalence $\Sigma \colon \T \longrightarrow \T$ and a class of triangles, called distinguished triangles, which are subject to the following axioms:

\textbf{(T1)} {\rm (i)} For any $X \in \T$, the triangle $\xymatrix{ X \ar[r]^{1} & X \ar[r] & 0 \ar[r] & \Sigma X}$ is distinguished.

{\rm (ii)} Any morphism $f \colon X \longrightarrow Y$ is part of a distinguished triangle $\xymatrix{ X \ar[r]^{f} & Y \ar[r]^{g} & Z \ar[r]^{h} & \Sigma X}$.

{\rm (iii)} Any triangle isomorphic to a distinguished one is itself distinguished.

\textbf{(T2)} A triangle $\xymatrix{ X \ar[r]^{f} & Y \ar[r]^{g} & Z \ar[r]^{h} & \Sigma X}$ is distinguished if and only if the triangle
$$\xymatrix{ Y \ar[r]^{g} & Z \ar[r]^{h} & \Sigma X \ar[r]^{-\Sigma f} & \Sigma Y}$$
is.

\textbf{(T3)} If the rows are distinguished triangles and the left square commutes in the diagram
$$\xymatrix{ X \ar[r]^f \ar[d]_{\alpha} & Y \ar[r]^g \ar[d]^{\beta} & Z \ar[r]^h \ar@{-->}[d]^{\gamma} & \Sigma X \ar[d]^{\Sigma \alpha} \\ X' \ar[r]^{f'} & Y' \ar[r]^{g'} & Z \ar[r]^{h'}  & \Sigma X',}$$
then there is a morphism $\gamma$ that makes the middle and right squares commute.

\textbf{(T4)}(Octahedral axiom) If $(f_1,g_1,h_1)$, $(f_2,g_2,h_2)$ and $(f_3,g_3,h_3)$ are distinguished triangles, $f_1$ and $f_2$ are composable and $f_3=f_2f_1$, then there exist morphisms $\alpha$ and $\beta$ such that $(\alpha, \beta, \Sigma g_1 \circ h_2)$ is a distinguished triangle and the diagram
$$\xymatrix{ X \ar[r]^{f_1} \ar@{=}[d] & Y \ar[r]^{g_1} \ar[d]_{f_2} & Z \ar[r]^{h_1} \ar@{-->}[d]^{\alpha} & \Sigma X \ar@{=}[d]\\ X \ar[r]^{f_3}  & W \ar[d]_{g_2} \ar[r]^{g_3}  & V \ar[r]^{h_3} \ar@{-->}[d]^{\beta} & \Sigma X \ar[d]^{\Sigma f_1} \\ & U \ar@{=}[r] \ar[d]_{h_2} & U \ar[d]^{\Sigma g_1 \circ h_2} \ar[r]^{h_2} & \Sigma Y \\ & \Sigma Y \ar[r]^{ \Sigma g_1} & \Sigma Z & }$$
commutes.

\end{defi}

Next we recall the notion of a homological functor.

\begin{defi} \label{def0.1} Let $\T$ be a triangulated category and $\A$ an abelian category. An additive functor $F \colon \T  \longrightarrow \A$ is said to be \emph{homological} if any distinguished triangle
$$\xymatrix{ X \ar[r]^{f} & Y \ar[r]^{g} & Z \ar[r]^{h} & \Sigma X}$$
yields a long exact sequence
$$\xymatrix{ \ldots \ar[r] & F(X) \ar[r]^{F(f)} & F(Y) \ar[r]^{F(g)} & F(Z) \ar[r]^{F(h)} & F(\Sigma X) \ar[r]^{F(\Sigma f)} & F(\Sigma Y) \ar[r] & \ldots .}$$
\end{defi}

As is shown in \cite[IV.1]{GM96}, the triangulated category axioms imply

\begin{prop} \label{ff} For any object $X$ of a triangulated category $\T$, the representable functor
$$\Hom_{\T}(X,-)_* \colon \T \longrightarrow \Grmod \Endd_{\T}(X)_*$$
is homological (Here $\Endd_{\T}(X)_*$ denotes the graded endomorphism ring of $X$.). \end{prop}

\begin{defi}\label{def0.3} Let $\T$ be a triangulated category with infinite coproducts. An object $P \in \T$ is \emph{compact} if the functor $\Hom_{\T}(P,-)$ commutes with arbitrary coproducts. One says that $P$ is a \emph{compact generator} if it is a compact object and in addition $\Hom_{\T}(P,X)_*=0$ implies that $X=0$, for any $X \in \T$ (or, equivalently, the functor $\Hom_{\T}(P,-)_* \colon \T \longrightarrow \Grmod \Endd_{\T}(P)_*$ reflects isomorphisms).  \end{defi}

The following lemma is well known (See e.g. \cite[Lemma 5.17 and Proposition A.4]{BKS04}. Although Lemma 5.17 in \cite{BKS04} is stated for the derived category of a differential graded
algebra, the proof given in \cite{BKS04} can be equally applied to general triangulated categories.).

\begin{lem} \label{proeqciulebi} Let $\T$ be a triangulated category with infinite coproducts and $P \in \T$ a compact generator.

{\rm (i)} If $F \in \T$ and $\Hom_{\T}(P,F)_*$ is a projective $\Endd_{\T}(P)_*$-module, then the map
$$ \Hom_{\T}(P,-)_* \colon\Hom_{\T}(F,X) \longrightarrow \Hom_{\Endd_{\T}(P)_*}(\Hom_{\T}(P,F)_*, \Hom_{\T}(P,X)_*)$$
is an isomorphism for any $X \in \T$.

{\rm (ii)} For any projective graded $\Endd_{\T}(P)_*$-module $M$, there exists $G \in \T$ such that
$$\Hom_{\T}(P,G)_* \cong M$$
in $\Grmod \Endd_{\T}(P)_*$.
\end{lem}

Next, using Lemma \ref{proeqciulebi}, one gets

\begin{prop}[The Adams spectral sequence] \label{ASS} Suppose $\T$ is a triangulated category with infinite coproducts and $P \in \T$ a compact generator, and assume that the graded global homological dimension of $\Endd_{\T}(P)_*$ is finite. Then there is a bounded convergent spectral sequence
$$E_2^{pq}=\Ext^p_{{\Endd_{\T}(P)}_*}(\Hom_{\T}(P,X)_*, \Hom_{\T}(P,Y)_*[q]) \Rightarrow \Hom_{\T}(X, \Sigma^{p+q}Y).$$
\end{prop}

In fact, Proposition \ref{ASS} is a special case of \cite[4.5]{C98}.

\vspace{0.75em}

\subsection{Model categories}~

\vspace{0.75em}

\noindent In this subsection we recall some basic properties of model categories and simplicial model categories.

\;

A \emph{model category} $\M$ is a bicomplete category equipped with three classes of morphisms called weak equivalences, fibrations and cofibrations, satisfying certain axioms. We will not list these axioms here. The point of this structure is that it allows one to ``do homotopy theory'' in $\M$. Good references for model categories include \cite{DS95} and \cite{H99}.

The fundamental example of a model category is the category \textbf{sSet} of simplicial sets (\cite{Q67}, \cite[1.11.3]{GJ99}). Further important examples are the category \textbf{Top} of topological spaces (\cite{Q67}, \cite[2.4.19]{H99}) and the category $\Ch(K)$ of chain complexes of modules over a ring $K$ \cite[2.3.11]{H99}.

For any model category $\M$, one has the associated \emph{homotopy category} $\Ho(\M)$ which is defined as the localization of $\M$ with respect to the class of weak equivalences (see e.g. \cite[1.2]{H99} or \cite{DS95}). The model structure guarantees that we do not face set theoretic problems when passing to localization, i.e., $\Ho(\M)$ has $\Hom$-sets. Note that $\Ho(\M)$ admits other equivalent descriptions as well. For example, Ho($\M$) is equivalent to the homotopy category of cofibrant objects. More precisely, let $\M_{\cof}$ denote the full subcategory of cofibrant objects in $\M$. Then Ho($\M$) is equivalent to the localization $\M_{\cof}[\W^{-1}]$, where $\W$ is the class of weak equivalences in $\M_{\cof}$. This description of the homotopy category is most convenient one for our exposition, and therefore we make the following
\begin{conv}\label{convI}In what follows, we let $\Ho(\M)$ denote the category $\M_{\cof}[\W^{-1}]$. \end{conv}
Note that the class $\W$ admits a calculus of left fractions \cite{Q67}. This is one of the advantages of Convention \ref{convI}.

\;

Further, recall the definition of a Quillen adjunction.
\begin{defi}\label{defI} A \emph{Quillen adjunction} between two model categories $\M$ and $\N$ is a pair of adjoint functors
$$\xymatrix{ F \colon \M \ar@<0.5ex>[r] & \N \cocolon G \ar@<0.5ex>[l]},$$
where the left adjoint $F$ preserves cofibrations and acyclic cofibrations (or, equivalently, $G$ preserves fibrations and acyclic fibrations).\end{defi}

We refer to $F$ as a left Quillen functor and to $G$ as a right Quillen functor. Quillen's total derived functor theorem (see e.g. \cite{Q67} or \cite[2.8.7]{GJ99}) says that any such pair of adjoint functors induces an adjunction
$$\xymatrix{ \textbf{L}F \colon \Ho(\M) \ar@<0.5ex>[r] & \Ho(\N) \cocolon \textbf{R}G \ar@<0.5ex>[l]}.$$
The functor $\textbf{L}F$ is called the left derived functor of $F$ and $\textbf{R}G$ the right derived functor of $G$. If $\textbf{L}F$  is an equivalence of categories (or, equivalently, $\textbf{R}G$ is an equivalence), then the Quillen adjunction is called a \emph{Quillen equivalence}.

Convention \ref{convI} allows us to give a very simple description of $\textbf{L}F$. Indeed, the functor
$$F\vert_{ \M_{\cof}} \colon \M_{\cof} \longrightarrow \N_{\cof}$$
preserves weak equivalences \cite[1.1.12]{H99} and, therefore, it induces a functor between the localizations. This functor is exactly $\textbf{L}F$ in terms of \ref{convI}.

An important class of model categories is the class of \emph{simplicial model categories} \cite[2.3]{GJ99}. These are model categories which are enriched, tensored and cotensored over \textbf{sSet} and which satisfy Quillen's axiom \textbf{SM7} \cite[2.3.1]{GJ99} (or, equivalently, the pushout product axiom \cite[2.3.11]{GJ99}). If a simplicial model category $\M$ is pointed, i.e., the terminal object is isomorphic to the initial one, then $\M$ is enriched over the category $\sSet _*$ of pointed simplicial sets.  In particular, we have the functors
$$- \wedge - \colon \sSet _* \times \M \longrightarrow \M, \;\; \;\; \Map_{\M}(-,-) \colon \M^{op} \times \M \longrightarrow \sSet_*.$$
and the adjunction
$$ \Hom_{\M}(K \wedge X,Y) \cong \Hom_{\sSet_*}(K, \Map_{\M}(X,Y)).$$

\begin{remk}\label{same letter} It follows from the axioms of pointed simplicial (model) categories that for any $n \geq 0$ and $X,Y \in \M$, there is a canonical natural isomorphism
$$\Map_{\M}(X,Y)_n \cong \Hom_{\M}(\Delta[n]_+ \wedge X, Y),$$
where $\Delta[n]_+$ is the union of $\Delta[n]$ with a disjoint base point. For any element $\theta$ in $\Map_{\M}(X,Y)_n$, the corresponding morphism $\Delta[n]_+ \wedge X \longrightarrow Y$ will be denoted again by $\theta$.

\end{remk}

\vspace{0.75em}

\subsection{Stable model categories}~

\vspace{0.75em}

\noindent Recall (\cite{Q67}, \cite[6.1.1]{H99}) that the homotopy category Ho($\M$) of a pointed model category $\M$ supports a \emph{suspension} functor
$$\Sigma \colon \Ho(\M) \longrightarrow \Ho(\M)$$
with a right adjoint \emph{loop} functor
$$\Omega \colon \Ho(\M) \longrightarrow \Ho(\M).$$

\begin{defi}\label{defII} A \emph{stable model category} is a pointed model category for which the functors $\Sigma$ and $\Omega$ are inverse equivalences.
\end{defi}

\begin{remk}\label{remkI}
If $\M$ is a pointed simplicial model category, then the suspension functor $$\xymatrix{\Sigma \colon \Ho(\M) \ar[r] & \Ho(\M)}$$ admits a very simple description. Indeed, by the simplicial model category axioms the functor
$$ S^1 \wedge - \colon \M \longrightarrow \M$$
is a left Quillen functor. Then $\Sigma$ can be defined as the left derived functor of $S^1 \wedge -$. \end{remk}

Note that if $\M$ is stable, then the homotopy category Ho($\M$) is a triangulated category \cite[7.1.6]{H99} with $\Sigma$ a shift functor. We will now recall the construction of distinguished triangles in Ho($\M$) in the case when $\M$ is simplicial.

\begin{defi}\label{cone} Suppose $\M$ is a simplicial stable model category and $f \colon X \longrightarrow Y$ a morphism in $\M_{\cof}$. The pushout of $f$ along the morphism $$ \incl \wedge 1 \colon X \cong S^0 \wedge X \longrightarrow (I,0) \wedge X$$ is called the \emph{mapping cone} of $f$ and is denoted by $\Cone(f)$. \end{defi}

Thus $\Cone(f)$ comes with the pushout square
$$\xymatrix{ X \ar[r]^f \ar[d] & Y \ar[d]^{\xi} \\ (I,0) \wedge X \ar[r]^{\mu} & \Cone(f).}$$
Furthermore, the universal property of pushout implies that there is a commutative diagram
$$\xymatrix{ Y \ar[r]^-{\xi} \ar[dr]_{\ast} & \Cone(f) \ar@{-->}[d]^{\delta} & (I,0) \wedge X \ar[l]_{\mu} \ar[dl]^{\pi \wedge 1} \\ &  S^1 \wedge X, & }$$
where $\pi \colon I \longrightarrow S^1$ is the projection.

\begin{defi}\label{defIV} Let $\M$ be a simplicial stable model category and $f \colon X \longrightarrow Y$ a morphism in $\M_{\cof}$. The \emph{elementary} triangle associated to $f$ is the triangle
$$\xymatrix{X \ar[r]^{f} & Y \ar[r]^-{\xi} & \Cone(f) \ar[r]^{\delta} & S^1 \wedge X.}$$
A triangle
$$\xymatrix{ A \ar[r] & B \ar[r] & C \ar[r] & \Sigma A}$$
in Ho($\M$) is called \emph{distinguished} if it is isomorphic to an elementary one.
\end{defi}

The distinguished triangles together with $\Sigma$ define the triangulated structure on Ho($\M$).

\;

Now let us give some examples of stable model categories which are important in what follows.

\begin{exmp}\label{exmpI} For any ring $K$, the category $\Ch(K)$ of unbounded chain complexes of $K$-modules with the projective model structure is a stable model category \cite[2.3.11]{H99}. The weak equivalences and fibrations in this model structure are quasi-isomorphisms (i.e., homology isomorphisms) and degreewise epimorphisms, respectively.
\end{exmp}

\begin{exmp}\label{exmpII} Let $R$ be a symmetric ring spectrum. The category $\Mod R$ of right $R$-modules admits a stable simplicial model structure \cite[Corollary 5.5.2]{HSS00}. The homotopy category Ho($\Mod R$) is called the derived category of $R$, denoted by $\D(R)$. Note that $R$ is a compact generator of $\D(R)$ and the representable homological functor associated to $R$ is the homotopy group functor $\pi_* \colon \D(R) \longrightarrow \Grmod \pi_*R$. \end{exmp}

\begin{exmp} \label{exmpIII} Let $A$ be a differential graded ring. Then, by Example \ref{exmpI} and \cite[4.1]{SS00}, the category $\Mod A$ of differential graded right modules over $A$ has a projective model structure in which the weak equivalences and fibrations are as in \ref{exmpI}. This model structure is stable as well. The category Ho($\Mod A$) is called the derived category of $A$, denoted by $\D(A)$. The differential graded module $A$ is a compact generator of $\D(A)$ and the representable homological functor associated to $A$ is the homology $H_* \colon \D(A) \longrightarrow \Grmod H_*A$.

By \cite[2.15]{S07}, there is a symmetric ring spectrum $HA$ such that $\Mod HA$ is Quillen equivalent to $\Mod A$. In particular, the derived category $\D(A)$ is triangulated equivalent to $\D(HA)$. Thus, from the point of view of model category theory Example \ref{exmpIII} is a special case of Example \ref{exmpII}. \end{exmp}

\begin{remk}\label{algcone} Note that the triangulated structure on the derived category $\D(A)$ of a differential graded ring $A$, comes from the usual shift functor and algebraic mapping cone construction. More precisely, for any $M \in \D(A)$, the suspension of $M$, denoted by $M[1]$, is given by
$$M[1]_n= M_{n-1}, \;\;\; d^{M[1]}_n= -d^{M}_{n-1}.$$
Obviously, $M[1]$ inherits the structure of a graded right $A$-module from that of $M$ (the action of $A$ on $M[1]$ is not twisted by a sign). Furthermore, it is clear that $d^{M[1]}$ satisfies the Leibniz rule.

Next, let $f \colon M \longrightarrow M'$ be a morphism of DG modules over $A$. The \emph{(algebraic) mapping cone} $C(f)$ of $f$ is defined by
$$(C(f))_n=M_{n-1} \oplus M'_n, \;\;\; n \in \Z,$$
$$d(m,m')=(-dm, f(m)+dm').$$
Clearly, $C(f)$ inherits the structure of graded (right) $A$-module from that of $M$ and $M'$, and the differential of $C(f)$ satisfies the Leibniz rule. Thus $(C(f), d)$ is a DG $A$-module. Further, there are canonical morphisms of DG modules
$$\iota \colon M'\longrightarrow C(f), \;\;\;\iota(m')=(0,m')$$
and
$$\partial \colon C(f) \longrightarrow M[1], \;\;\; \partial(m,m')=m.$$
These morphisms together with $f \colon M \longrightarrow M'$ form a triangle
$$\xymatrix{ M \ar[r]^f & M' \ar[r]^\iota & C(f) \ar[r]^\partial & M[1],}$$
called the elementary triangle associated to $f$.

A triangle in $\D(A)$ is distinguished if and only if it is isomorphic to an elementary one.

\end{remk}

\vspace{0.75em}

\subsection{Diagram categories}~

\vspace{0.75em}

\noindent We will use model structures on diagram categories throughout the paper. Therefore, we recall some facts and notions concerning them.

\begin{defi}[{\cite[5.1.1]{H99}}]\label{defV} Let $\omega$ denote the poset category of the ordered set  $\{0,1,2,...\}$. A small category $\C$ is called a \emph{direct category} if there exists a functor $f\colon \C \longrightarrow \omega$ that sends non-identity morphisms to non-identity morphisms. \end{defi}

\begin{defi}[{\cite[5.1.2]{H99}}]\label{defVI}Suppose $\M$ is a category with small colimits, $\C$ a small category, $z$ an object in $\C$ and $\C_z$ the category of all non-identity morphisms in $\C$  with codomain $z$. The \emph{latching space functor} $L_z \colon \M^{\C} \longrightarrow \M$ is the composition
$$\xymatrix{ \M^{\C} \ar[r] & \M^{\C_z} \ar[r]^{\colim} & \M,}$$
where the first arrow is restriction. Note that we have a natural transformation $$ L_zX \longrightarrow X_z$$ for any fixed object $z \in \C$.
\end{defi}

The following proposition is proved in \cite[5.1.3]{H99}.

\begin{prop}\label{Reedy} Given a model category $\M$ and a direct category $\C$, there is a model structure on $\M^{\C}$ in which a morphism of diagrams $f\colon X \longrightarrow Y$ is a weak equivalence (resp. fibration) if and only if the map $f_z \colon X_z \longrightarrow Y_z$ is a weak equivalence (resp. fibration) for all $z \in \C$. Furthermore, $f\colon X \longrightarrow Y$ is an (acyclic) cofibration if and only if the induced map $X_z \amalg_{L_zX} L_zY \longrightarrow Y_z$ is an (acyclic) cofibration for all $z \in \C$. \end{prop}

\begin{exmp}\label{exmpIV} Any finite poset $\P$ is a direct category. Therefore, for any model category $\M$, we have the model structure \ref{Reedy} on $\M^{\P}$. \end{exmp}

\begin{exmp}\label{exmpV} Suppose $\M$ is a model category and $N \geq 2$ a natural number. Let $\C_N$ denote the poset consisting of elements $\{\beta_i$, $\zeta_i | i \in \Z/N\Z  \}$ such that $\beta_i < \zeta_i$ and $\beta_i < \zeta_{i+1}$, i.e.,
$$\xymatrix{\zeta_0
  \ar@{<-}[drrrr]%
    |<<<<<<<<<<<<{\text{\Large \textcolor{white}{$\blacksquare$}}}%
    |<<<<<<<<<<<<<<<<<<<<<{\text{\Large \textcolor{white}{$\blacksquare$}}}
    |>>>>>>>>>>>>>>>>>>>>>>>{\text{\Large \textcolor{white}{$\blacksquare$}}}
    |>>>>>>>>>>>>>>{\text{\Large \textcolor{white}{$\blacksquare$}}}
    |>>>>>>>>>>{\text{\Large \textcolor{white}{$\blacksquare$}}}%
& \zeta_1 & \ldots & \zeta_{N-2} & \zeta_{N-1} \\ \beta_0 \ar[u] \ar[ur] & \beta_1 \ar[u] \ar[ur] & \ldots \ar[ur] & \beta_{N-2} \ar[u] \ar[ur] & \beta_{N-1}, \ar[u] }$$
and consider $X \in \M^{\C_N}$. Then, by Proposition \ref{Reedy}, $X$ is cofibrant if and only if the canonical map $$ L_zX \longrightarrow X_z$$ is a cofibration in $ \M$ for all $z \in \C_N$, i.e., if and only if $X_{\beta_i}$ is cofibrant and the induced morphism $X_{\beta_{i-1}} \vee X_{\beta_i} \longrightarrow X_{\zeta_i}$ is a cofibration, for any $i \in \Z/N\Z$. \end{exmp}

\begin{exmp}\label{exmpVI} The category
$$\xymatrix{ \bullet \ar@<0.5ex>[r] \ar@<-0.5ex>[r] & \bullet} $$
is a direct category. Thus, by Proposition \ref{Reedy}, we get a model structure on the category of parallel arrows in $\M$.
\end{exmp}

It follows from \cite[5.1.5]{H99} that for any model category $\M$ and direct category $\C$, there is a Quillen adjunction
$$\xymatrix{ \colim\colon \M^{\C} \ar@<0.5ex>[r] & \M \cocolon \const, \ar@<0.5ex>[l],}$$
where $\colim$ is a left Quillen functor.

\begin{defi}\label{hocolim} The left derived functor of $\colim \colon \M^{\C} \longrightarrow \M$ is called the homotopy colimit, denoted by
$$\Hocolim \colon \Ho(\M^{\C}) \longrightarrow \Ho(\M).$$
If $\C$ is as in Example \ref{exmpVI}, then the homotopy colimit is called the homotopy coequalizer.
\end{defi}

\begin{remk}\label{diagmap} Let $\C$ be a direct category (Definition \ref{defV}). If $\M$ is a simplicial model category, then so is $\M^{\C}$. Indeed, the mapping spaces for $\M^{\C}$ are given by the end construction \cite[IX.5]{M98}

$$\Map_{\M^{\C}}(X,Y)=\int_{z \in \C} \Map_{\M}(X_z, Y_z).$$ The tensors and cotensors are defined levelwise. \end{remk}

\vspace{0.75em}

\subsection{Long exact sequences}~

\vspace{0.75em}

\noindent In this subsection we recall the Mayer-Vietoris sequence of a homotopy pullback and the long exact sequence of a homotopy coequalizer.

\begin{lem}[\cite{DR80}]\label{MV} Let
$$\xymatrix{ A \ar[r]^{\alpha} \ar[d]_{\beta} & B \ar[d]^{\phi} \\ C \ar[r]^{\psi} & D}$$
be a homotopy pullback square of pointed simplicial sets. Then there is a long exact Mayer-Vietoris sequence
\begin{multline*} \ldots \xto{}  \pi_n(A) \xto{(\alpha_*, \beta_*)} \pi_n(B)  \times \pi_n(C) \xto{(\phi_*, -\psi_*)} \pi_n(D) \xto{\partial} \pi_{n-1}(A) \xto{} \ldots \\ \ldots \xto{} \pi_1(A) \xto{(\alpha_*, \beta_*)} \pi_1(B)  \times \pi_1(C) \xto{(\phi_*, \psi_*^{-1})}  \pi_1(D) \xto{\partial} \pi_0(A) \xto{(\alpha_*, \beta_*)} \pi_0(B)  \times \pi_0(C).  \end{multline*} \end{lem}

\;

\;

Note that if
$$\xymatrix{ A \ar[r]^{\alpha} \ar[d]_{\beta} & B \ar[d]^{\phi} \\ C \ar[r]^{\psi} & D}$$
is a pullback square of Kan simplicial sets and $\phi$ a Kan fibration, then it is in fact a homotopy pullback. In this case the connecting homomorphism
$$\partial \colon \pi_n(D) \longrightarrow \pi_{n-1}(A)$$
can be described as follows.

\begin{lem}\label{dMV} Suppose that
$$\xymatrix{ A \ar[r]^{\alpha} \ar[d]_{\beta} & B \ar[d]^{\phi} \\ C \ar[r]^{\psi} & D}$$
is a pullback square of pointed Kan simplicial sets, and $\phi$ a Kan fibration. Suppose further that $F$ is the fiber of $\phi$ and $\beta$, and $\varepsilon \colon F \longrightarrow A$ the inclusion. Then the connecting homomorphism
$$\partial \colon \pi_n(D) \longrightarrow \pi_{n-1}(A)$$
of the Mayer-Vietoris sequence is given by the composite
$$\xymatrix{ \pi_n(D) \ar[r]^{\partial'} & \pi_{n-1}(F) \ar[r]^{\varepsilon_*} & \pi_{n-1}(A),}$$
where $\partial'$ is the connecting homomorphism of the long exact sequence of the fibration $\phi\colon B \longrightarrow D$  \cite[3.5]{K58}. \end{lem}

\begin{lem}\label{MVHspace} Let
$$\xymatrix{ A \ar[r]^{\alpha} \ar[d]_{\beta} & B \ar[d]^{\phi} \\ C \ar[r]^{\psi} & D}$$
be a homotopy pullback square of simplicial sets with $A$, $B$, $C$, $D$ homotopy associative and homotopy commutative $H$-spaces, and $\alpha$, $\beta$, $\phi$, $\psi$ $H$-space maps. Then the long exact Mayer-Vietoris sequence \ref{MV} can be extended to the right by one further term. More precisely, there is a long exact sequence of abelian groups

\begin{multline*} \ldots \xto{}  \pi_n(A) \xto{(\alpha_*, \beta_*)} \pi_n(B)  \oplus \pi_n(C) \xto{(\phi_*, -\psi_*)} \pi_n(D) \xto{\partial} \pi_{n-1}(A) \xto{} \ldots \\ \ldots \xto{} \pi_1(B)  \oplus \pi_1(C) \xto{(\phi_*, -\psi_*)}  \pi_1(D) \xto{\partial} \pi_0(A) \xto{(\alpha_*, \beta_*)} \pi_0(B)  \oplus \pi_0(C) \xto{(\phi_*, -\psi_*)} \pi_0(D).  \end{multline*} \end{lem}

Next we recall the following

\begin{lem}\label{htpycoeqtr} Let $\M$ be a stable model category and
$$\xymatrix{ X \ar@<0.5ex>[r]^f \ar@<-0.5ex>[r]_g & Y \ar[r]^h & Z}$$
a coequalizer diagram in $\M$. Suppose that $X$ is cofibrant and
$$ (f,g) \colon X \vee X \lra Y$$
is a cofibration (i.e., the given coequalizer is a homotopy coequalizer). Then there is a morphism $\delta\colon Z \longrightarrow \Sigma X$ in $\Ho(\M)$ such that the triangle
$$\xymatrix{ X \ar[r]^{f-g} & Y \ar[r]^h & Z \ar[r]^\delta & \Sigma X}$$
is distinguished.
\end{lem}

This lemma essentially follows from the remark at the end of Section I.4 in \cite{B73}.

\begin{coro}\label{htpycoeq} Let $\M$ be a stable model category, $\A$ an abelian category and $F \colon \Ho(\M)  \longrightarrow \A$ a homological functor. Then any homotopy coequalizer
$$\xymatrix{ X \ar@<0.5ex>[r]^f \ar@<-0.5ex>[r]_g & Y \ar[r]^h & Z}$$
in $\M$ induces a long exact sequence
$$\ldots \xto{} F(X) \xto{F(f)-F(g)} F(Y) \xto{F(h)} F(Z) \xto{\partial} F(\Sigma X) \xto{} \ldots.$$ \end{coro}

\section{Franke's functor}~\label{yv. mtavari}

\subsection{Formulation of main results}~

\vspace{0.75em}

\noindent Suppose $\M$ is a stable model category and $P \in \Ho(\M)$ a compact generator. In the sequel we use the notation $$\pi_*X = [P,X]_*=\Hom_{\Ho(\M)}(\Sigma^*P,Y).$$ In particular, $\pi_*P$ denotes the graded ring $[P,P]_*$.

\begin{remk} \label{wanacvleba} The suspension functor $\Sigma \colon \Ho(\M) \longrightarrow \Ho(\M)$ induces a natural isomorphism of graded $\pi_*P$-modules
$$\pi_*(\Sigma X) \cong (\pi_*X)[1].$$
Having this in mind, we identify $\pi_*(\Sigma X)$ with $(\pi_*X)[1]$ in what follows.
\end{remk}

\begin{defi} Suppose $N \geq 2$ is a natural number. A graded ring $B$ is said to be $N$-sparse if it is concentrated in degrees divisible by $N$. \end{defi}

\begin{defi} Let $\M$ be a stable model category and $P \in \M$ a compact generator. We say that an object $X \in \Ho(\M)$ is $k$-dimensional (with respect to $P$) if the projective dimension of $\pi_*X$ in $\Grmod \pi_*P$ is equal to $k$.
\end{defi}

Now we are ready to list the main results and their corollaries.

\begin{theo}\label{I}
Let $\M$ be a stable model category and $P$ a fixed compact generator of $\Ho(\M)$. Suppose that the graded ring $\pi_*P$ is $N$-sparse and the (right) graded global homological dimension of $\pi_*P$ is less than $N$. Then one can construct a functor
$$\R \colon \D(\pi_*P) \longrightarrow \Ho(\M)$$
such that the following hold:

{\rm (i)} The diagram $$\xymatrix{\D(\pi_*P) \ar[rr]^-{\R} \ar[dr]_-{H_*} & & \Ho(\M) \ar[dl]^-{\pi_*} \\ & \Grmod{\pi_*P} & }$$ commutes up to a natural isomorphism.

{\rm (ii)} The functor $\R$ restricts to an equivalence of the full subcategories of at most one-dimensional objects.

{\rm (iii)} Let $X \in \Ho(\M)$ and suppose that $X$ has dimension at most two. Then there exists $M \in \D(\pi_*P)$ such that $\R(M)$ is isomorphic to $X$ in $\Ho(\M)$. \end{theo}

\begin{theo}\label{equival} Let $\M$ be a stable model category and $P$ a compact generator of $\Ho(\M)$. Suppose that the graded ring $\pi_*P$ is $N$-sparse and the (right) graded global homological dimension of $\pi_*P$ is less than $N-1$. Then the functor $$\R \colon \D(\pi_*P) \longrightarrow \Ho(\M)$$ restricts to an equivalence of the full subcategories of at most two dimensional objects. \end{theo}

In particular, we have

\begin{coro}Let $\M$ be a stable model category and $P$ a compact generator of $\Ho(\M)$. Suppose that the graded ring $\pi_*P$ is $N$-sparse, the (right) graded global homological dimension of $\pi_*P$ is equal to two and $N \geq 4$. Then the functor $$\R \colon \D(\pi_*P) \longrightarrow \Ho(\M)$$ is an equivalence of categories.

\end{coro}

Theorems  \ref{Rcoro} and \ref{equival1} are obtained as special cases of Theorem \ref{I} and Theorem \ref{equival}, respectively. Further, in view of \ref{exmpIII}, one has

\begin{coro}\label{dgcoro} Let $A$ be a differential graded ring. Suppose that the graded homology ring $H_*A$ is concentrated in dimensions divisible by a natural number $N$ and assume that the (right) graded global homological dimension of $H_*A$ is less than $N$.
Then one can construct a functor
$$\R \colon \D(H_*A) \longrightarrow \D(A)$$
such that the following hold:

{\rm (i)} The diagram
$$\xymatrix{\D(H_*A) \ar[rr]^-{\R} \ar[dr]_{H_*} & & \D(A) \ar[dl]^-{H_*} \\ & \Grmod{H_*A} & }$$
commutes up to a natural isomorphism.

{\rm (ii)} The functor $\R \colon \D(H_*A) \longrightarrow \D(A)$ restricts to an equivalence of the full subcategories of at most one-dimensional objects.

{\rm (iii)} Let $X \in \D(A)$ and suppose $X$ has dimension at most two. Then there exists $M \in \D(H_*A)$ such that $\R(M)$ is isomorphic
to $X$ in $\D(A)$.

\end{coro}

\begin{coro}\label{dgequival1} Let $A$ be a differential graded ring. Suppose that the graded homology ring $H_*A$ is $N$-sparse and the (right) graded global homological dimension of $H_*A$ is less than $N-1$. Then the functor $$\R \colon \D(H_*A) \longrightarrow \D(A)$$ restricts to an equivalence of the full subcategories of at most two-dimensional objects. \end{coro}

\begin{coro}\label{dgequival2} Let $A$ be a differential graded ring. Suppose that the graded homology ring $H_*A$ is $N$-sparse, the (right) graded global homological dimension of $H_*A$ is equal to two and $N \geq 4$. Then the functor $$\R \colon \D(H_*A) \longrightarrow \D(A)$$  is an equivalence of categories.
\end{coro}

The functor $\R$ is constructed in \ref{Q} and the proof of \ref{I} (i) is given in \ref{R}. Next, \ref{I} (ii) is proved in Section \ref{advgamoy}, and finally, the proofs of \ref{I} (iii) and \ref{equival} are discussed in Section \ref{intgamoy}.

\vspace{0.75em}

\subsection{Some consequences of the Adams spectral sequence}~\label{spectral}

\vspace{0.75em}

\noindent In this subsection we state some known facts which follow from the Adams spectral sequence \ref{ASS}
$$E_2^{pq}=\Ext^p_{\pi_*P}(\pi_*X, \pi_*Y[q]) \Rightarrow [X, \Sigma^{p+q}Y].$$

Let $\M$, $P$ and $N$ be as in \ref{I}. Then the category of graded $\pi_*P$-modules splits as follows

$$\Grmod{\pi_*P} \sim \B \oplus \B[1] \oplus \ldots \oplus \B[N-1],$$ where $\B$ is the full subcategory of $\Grmod{\pi_*P}$ consisting of all those modules which are concentrated in degrees divisible by $N$.

Let $\E_i$ denote the full subcategory of $\Ho(\M)$ consisting of objects $X \in \Ho(\M)$ with $$\pi_*X \in \B[i].$$

\begin{prop} \label{II}

{\rm (i)} The functor
$$\pi_*\vert_{ \E_i} \colon \E_i \longrightarrow \B[i]$$
is an equivalence of categories. In particular, for any $i \in \Z/N\Z$ and $X, Y \in \E_i$, the natural map $$\pi_* \colon [X,Y] \longrightarrow \Hom_{\pi_*P}(\pi_*X, \pi_*Y)$$ is an isomorphism.

{\rm (ii)} For any $i \in \Z/N\Z$, $X \in \E_{i-1}$ and $Y \in \E_i$, there is a natural isomorphism $$[X, Y] \cong \Ext^{1}_{\pi_*P}(\pi_*X[1], \pi_*Y).$$

{\rm (iii)} For any $i \in \Z/N\Z$ and $X, Y \in \E_i$, there is a natural isomorphism $$[\Sigma X, Y] \cong \Ext^{N-1}_{\pi_*P}(\pi_*X[N], \pi_*Y).$$ \end{prop}

This proposition immediately follows from the Adams spectral sequence. (The conditions in \ref{I} imply that the corresponding spectral sequences collapse at $E_2$ and all appropriate terms vanish.)

\begin{remk}\label{III} If $X \in \E_{i-1}$ and $Y \in \E_i$ then $\pi_*f=0$ for any $f \in [X,Y]$. This and the construction of the Adams spectral sequence enables one to get an explicit description of the isomorphism in \ref{II}(ii). Indeed, for any $f \in [X,Y]$, a distinguished triangle
$$\xymatrix{X \ar[r]^-{f} & Y \ar[r]^-{\xi} & \Cone (k) \ar[r]^-{\delta} & \Sigma X}$$
gives a short exact sequence
$$\xymatrix{E \colon 0 \ar[r] & \pi_*Y  \ar[r]^-{\pi_*\xi} & \pi_*\Cone (k) \ar[r]^{\pi_*\delta} & \pi_*X[1] \ar[r] & 0,}$$
and the isomorphism \ref{II}(ii) sends $f$ to the extension class of $E$.
\end{remk}

\vspace{0.75em}

\subsection{Crowned diagrams and DG modules}~\label{Q}

\vspace{0.75em}

\noindent In this subsection we start with the proof of Theorem \ref{I}. In order to simplify the exposition, the model category $\M$ will be assumed to be simplicial in the sequel. (Note that our examples of model categories (see \ref{Rcoro} and \ref{dgcoro}) are in fact simplicial or Quillen equivalent to simplicial ones.) However, the proof can be applied for stable model categories without any enrichment as well, using the techniques of cosimplicial frames from \cite{H99} (see \ref{csframes} for details).

Let $\C_N$ be as in \ref{exmpV}. Suppose $X$ is an object of $(\M^{\C_N})_{\cof}$ and
$$l_i \colon X_{\beta_i} \longrightarrow X_{\zeta_i}, \;\;\;\;\;\;\;\;\;\;\;\;\; k_i \colon X_{\beta_{i-1}} \longrightarrow X_{\zeta_i},\;\;\;\;\;\;\;\;\; i \in \Z/N\Z $$
the morphisms of $X$. Since $X$ is cofibrant in $\M^{\C_N}$, the objects $X_{\beta_i}$ and $X_{\zeta_i}$ are cofibrant in $\M$, $i \in \Z/N\Z$ . Let
$$ Z^{(i)}(X)= \pi_*X_{\zeta_i},\;\;\;\;\;\;\;\; B^{(i)}(X)= \pi_*X_{\beta_i}\;\;\;\;\;\;\;\; C^{(i)}(X)= \pi_*\Cone(k_i),$$ $$\lambda^{(i)}=\pi_*l_i \colon B^{(i)}(X) \longrightarrow Z^{(i)}(X),\;\;\;\;\;\;\;\;i \in \Z/N\Z,$$
where $\Cone(k_i)$ denotes the cone construction from \ref{cone} for simplicial model categories applied to $k_i$. Consider the mapping cone sequence
$$\xymatrix{X_{\beta_{i-1}} \ar[r]^{k_i} & X_{\zeta_i} \ar[r] & \Cone(k_i) \ar[r] & \Sigma X_{\beta_{i-1}}}$$
of $k_i$, $i \in \Z/N\Z$. Applying $\pi_* \colon \Ho(\M) \longrightarrow \Grmod \pi_*P$ to
$\xymatrix{X_{\zeta_i} \ar[r] & \Cone(k_i) \ar[r] & \Sigma X_{\beta_{i-1}},}$
we obtain the sequence $$ Z^{(i)}(X) \xto{\iota^{(i)}} C^{(i)}(X) \xto{\rho^{(i)}} B^{(i-1)}(X)[1]$$ (see \ref{wanacvleba}). Clearly, $Z^{(i)}$ and $B^{(i)}$ are functors from $(\M^{\C_N})_{\cof}$ to $\Grmod \pi_*P$ and $\lambda^{(i)}$ is natural. Besides, since the cone construction of \ref{cone} is functorial, $C^{(i)}$ is a functor from $(\M^{\C_N})_{\cof}$ to $\Grmod \pi_*P$ as well, and $\iota^{(i)}$ and $\rho^{(i)}$ are natural. Further, $Z^{(i)}$, $B^{(i)}$ and $C^{(i)}$ pass through Ho($\M^{\C_N}$) (as they send weak equivalences of diagrams to isomorphisms), and $\lambda^{(i)}$, $\iota^{(i)}$, $\rho^{(i)}$ induce natural transformations on the homotopy level. We use the same notation for the resulting functors and natural transformations.

Next, let $\L$ be the full subcategory of Ho($\M^{\C_N}$) consisting of all those $X$ which satisfy the following conditions:

\;

{\rm (i)} The graded $\pi_*P$-modules $Z^{(i)}(X)$, $B^{(i)}(X)$ are objects of $\B[i]$ for any $\in \Z/N\Z$;

\;

{\rm (ii)} The morphism $\lambda^{(i)}$ is a monomorphism for all $i$ $\in$ $\Z/N\Z$.

We construct a functor $$\Q \colon \L \longrightarrow \Mod \pi_*P.$$

Let $X$ be an object of $\L$. As the functor
$$ \pi_* \colon \Ho(\M) \longrightarrow \Grmod{\pi_*P}$$
is homological, the distinguished triangles
$$\xymatrix{X_{\beta_{i-1}} \ar[r]^{k_i} & X_{\zeta_i} \ar[r] & \Cone(k_i) \ar[r] & \Sigma X_{\beta_{i-1}}}$$
induce long exact sequences
$$ \ldots \longrightarrow B^{(i-1)}(X) \longrightarrow Z^{(i)}(X)  \xto{\iota^{(i)}} C^{(i)}(X) \xto{\rho^{(i)}} B^{(i-1)}(X)[1] \longrightarrow Z^{(i)}(X)[1] \longrightarrow \ldots. $$
Note that $B^{(i-1)}(X)$ $\in$ $\B[i-1]$ and $Z^{(i)}(X)$ $\in$ $\B[i]$ for all $i$ $\in$ $\Z/N\Z$, since $X$ $\in$ $\L$. Therefore, the morphisms $B^{(i-1)}(X) \longrightarrow Z^{(i)}(X)$ and $B^{(i-1)}(X)[1] \longrightarrow Z^{(i)}(X)[1]$ are zero. Consequently, for any $i \in \Z/N\Z$,  we get a short exact sequence $$\xymatrix{ 0 \ar[r] & Z^{(i)}(X)  \ar[r]^{\iota^{(i)}} & C^{(i)}(X) \ar[r]^{\rho^{(i)}} & B^{(i-1)}(X)[1] \ar[r] & 0}$$ in $\Grmod{\pi_*P}$.  Then, consider the following graded $\pi_*P$-modules
$$C_*(X)= C^{(0)}(X) \oplus C^{(1)}(X) \oplus \ldots \oplus C^{(N-1)}(X),$$ $$Z_*(X)= Z^{(0)}(X) \oplus Z^{(1)}(X) \oplus \ldots \oplus Z^{(N-1)}(X),$$ $$B_*(X)= B^{(0)}(X) \oplus B^{(1)}(X) \oplus \ldots \oplus B^{(N-1)}(X).$$
The morphisms $\lambda^{(i)}$, $\iota^{(i)}$, $\rho^{(i)}$, $i$ $\in$ $\Z/N\Z$, induce morphisms between the direct sums
\begin{align*}
\lambda &: B_*(X) \longrightarrow Z_*(X), & \lambda &= \lambda^{(0)} \oplus \lambda^{(1)} \oplus \ldots \oplus \lambda^{(N-1)}, \\
  \iota &: Z_*(X) \longrightarrow C_*(X), &   \iota &=   \iota^{(0)} \oplus   \iota^{(1)} \oplus \ldots \oplus   \iota^{(N-1)}, \\
   \rho &: C_*(X) \longrightarrow B_*(X)[1], & \rho &=    \rho^{(0)} \oplus    \rho^{(1)} \oplus \ldots \oplus \rho^{(N-1)}.
\end{align*}
After summing up, we get a short exact sequence of $\pi_*P$-modules
\begin{equation}\label{se}\begin{split}\xymatrix{0 \ar[r] & Z_*(X)   \ar[r]^{\iota} & C_*(X) \ar[r]^{\rho} & B_*(X)[1] \ar[r] & 0.}\end{split}
\end{equation}
Splicing this short exact sequence with its shifted copy gives a differential graded $\pi_*P$-module. More precisely, define $$d=\iota[1] \lambda[1] \rho \colon C_*(X) \longrightarrow C_*(X)[1].$$ Then $dd=0$ and, therefore, we get a DG $\pi_*P$-module $(C_*(X), d)$. The desired functor $$\Q \colon \L \longrightarrow \Mod \pi_*P,$$ is defined by $\Q(X)=(C_*(X), d)$.

\begin{remk}\label{IV} Let $\L'$ denote the full subcategory of $(\M^{\C_N})_{\cof}$ consisting of all the objects $X$ with $Z^{(i)}(X)$, $B^{(i)}(X) \in$ $\B[i]$ and $\lambda^{(i)}$ a monomorphism, for any $i \in \Z/N\Z$ (i.e., $\L'$ has the same objects as $\L$). Then $\L$ is the localization of $\L'$ with respect to the class of weak equivalences. This immediately follows from the fact that the class of weak equivalences in $(\M^{\C_N})_{\cof}$ admits a calculus of left fractions \cite{Q67}. \end{remk}

\begin{prop}\label{V} Let $\K$ be the full subcategory of $\L$ consisting of all the objects $X$ with $$\proj  Z^{(i)}(X)<N-1,\;\;\;\;\;\;\;\;\;\;\;\;\;\;\;\ \proj  B^{(i)}(X)<N-1$$
for any $i \in\Z/N\Z$ (see \ref{grdims}). Then the functor $$\Q\vert_{ \K} \colon \K \longrightarrow \Mod \pi_*P$$ is full and faithful.\end{prop}

\begin{prf} Take $X$, $\widetilde{X} \in \K$. Let $K$ denote the abelian group
$$\Hom_{\Ho(\M^{\C_N})}(X,\widetilde{X})=[X, \wt{X}],$$
and $L$ the abelian group of commutative diagrams of the form
$$\xymatrix{Z_*(X) \ar[r]^-\iota \ar[d] & C_*(X) \ar[d] \\ Z_*(\widetilde{X}) \ar[r]^-{\widetilde{\iota}} & C_*(\widetilde{X}).}$$
In other words,
$$L=\Hom_{(\Grmod \pi_*P)^{\underline{1}}}(Z_*(X)  \xto{\iota} C_*(X), Z_*(\widetilde{X})  \xto{\widetilde{\iota}} C_*(\widetilde{X})),$$
where $\underline{1}$ denotes the category $\bullet \longrightarrow \bullet$. To prove the proposition it suffices to check that the morphism $$q \colon K \longrightarrow L$$ induced by the functor $\Q$ is injective and its image consists of all those morphisms which respect the differentials, i.e., which are morphisms of DG modules. Thanks to (\ref{se}) any $f \in L$ induces a morphism on the cokernels $$B_*(X) \xto{\overline{f}} B_*(\widetilde{X}).$$ By definition of $d \colon C_*(X) \longrightarrow C_*(X)[1]$, a morphism $ f$ $\in$ $L$ is a DG morphism if and only if the outer square in the diagram
$$\xymatrix{C_*(X) \ar[r]^{\rho} \ar[d]^{f} & B_*(X)[1] \ar[r]^{\lambda[1]} \ar[d]^{\ol{f}[1]} & Z_*(X)[1] \ar[r]^{\iota[1]} \ar[d]^{f[1]} & C_*(X)[1] \ar[d]^{f[1]} \\ C_*(\widetilde{X}) \ar[r]^{\widetilde{\rho}} & B_*(\widetilde{X})[1] \ar[r]^{\widetilde{\lambda}[1]} & Z_*(\widetilde{X})[1] \ar[r]^{\widetilde{\iota}[1]} & C_*(\widetilde{X})[1]}$$
commutes. The left and right squares in the diagram are commutative. Therefore, since $\wt{\iota}[1] \colon Z_*(\widetilde{X})[1] \longrightarrow C_*(\widetilde{X})[1]$ is a monomorphism and $\rho \colon C_*(X) \longrightarrow B_*(X)[1]$ an epimorphism, the outer square commutes if and only if the middle one does, i.e., if and only if $f$ lies in the kernel of
$$ D \colon L \longrightarrow \Hom_{\pi_*P}(B_*(X), Z_*(\widetilde{X})),\;\;\;\;\; {D(f)}=\widetilde{\lambda}  {\overline{f}} - {f \lambda}.$$
Thus, to show that the functor $$\Q\vert_{ \K} \colon \K \longrightarrow \Mod \pi_*P$$ is full and faithful is equivalent to showing that the sequence
\begin{align} \begin{split}\label{eseq}\xymatrix{ 0 \ar[r] & K \ar[r]^{q} & L \ar[r]^-{D} & \Hom_{\pi_*P}(B_*(X), Z_*(\widetilde{X}))}\end{split} \end{align} is exact. To prove the latter we describe $K$ and $L$ in terms of some exact sequences.

Let us start with $K$. By \ref{diagmap} we have a pullback square of based simplicial sets
\begin{align} \begin{split}\label{plbk} \xymatrix{\Map_{\M^{\C_N}}(X,\widetilde{X})  \ar[rrr] \ar[d]  & & &   \us{i }{\prod} \Map_{\M}(X_{\zeta_i}, \widetilde{X}_{\zeta_i}) \ar[d]^{\phi} \\ \us{i}{\prod} \Map_{\M}(X_{\beta_i}, \widetilde{X}_{\beta_i}) \ar[rrr]^-{\psi} & & &  \us{i}{\prod} \Map_\M(X_{\beta_{i-1}}, \widetilde{X}_{\zeta_i}) \times \us{i}{\prod} \Map_\M(X_{\beta_i}, \widetilde{X}_{\zeta_i}), }\end{split} \end{align}
where ${\phi}$ and $\psi$ are induced by
$$X_{\beta_{i-1}} \vee X_{\beta_{i}} \xto{(k_i,\; l_i)} X_{\zeta_i}\;\; \; \; \; \;  \text{and} \;\; \; \; \; \;  \widetilde{X}_{\beta_i} \xto{(\wt{k}_{i+1},\; \wt{l}_i)} \widetilde{X}_{\zeta_{i+1}} \times \widetilde{X}_{\zeta_i},  \; \; i \in \Z/N\Z,$$
respectively. Without loss of generality we may assume that $\widetilde{X}$ is fibrant in $\M^{\C_N}$, i.e., $\widetilde{X}_{\beta_i}$ and $\widetilde{X}_{\zeta_i}$, $i \in \Z/N\Z$ are fibrant in $\M$. Besides, since $X$ is a cofibrant object in $\M^{\C_N}$, we can conclude that
$$ X_{\beta_{i-1}} \vee X_{\beta_{i}} \xto{(k_i,\; l_i)} X_{\zeta_i}$$
is a cofibration and $X_{\beta_{i}}$, $X_{\zeta_i}$ are cofibrant, for any $i \in \Z/N\Z$ (see \ref{exmpV}). This, by the axiom \textbf{SM7} for simplicial model categories implies that all objects in (\ref{plbk}) are fibrant and $\phi$ is a Kan fibration. Next, as $\M$ is stable, the simplicial sets and maps involved in (\ref{plbk}) are infinite loop spaces and infinite loop maps, respectively. In particular, all the objects are homotopy associative and homotopy commutative $H$-spaces and all the morphisms are $H$-space maps. Thus, by \ref{MVHspace}, the diagram (\ref{plbk}) gives a long exact Mayer-Vietoris sequence. Let us identify the terms of this long exact sequence. Since $X$ is cofibrant and $\widetilde{X}$ is fibrant,
the group $\pi_0(\Map_{\M^{\C_N}}(X,\widetilde{X}))$ is isomorphic to $[X, \widetilde{X}]=K$, by \cite[2.3.10]{GJ99}. As noted above, $X_{\beta_i}$ $X_{\zeta_i}$, $i$ $\in$ $\Z/N\Z$ are cofibrant and $\widetilde{X}_{\beta_i}$, $\widetilde{X}_{\zeta_i}$ are fibrant for all $i \in \Z/N\Z$. Hence, using again \cite[2.3.10]{GJ99}, we get isomorphisms
\begin{align*}
\pi_n(\Map_{\M}(X_{\beta_i}, \widetilde{X}_{\beta_i})) &\cong [\Sigma^n X_{\beta_i}, \widetilde{X}_{\beta_i}],
\\
\pi_n(\Map_{\M}(X_{\zeta_i}, \widetilde{X}_{\zeta_i})) &\cong [\Sigma^n X_{\zeta_i}, \widetilde{X}_{\zeta_i}],
\\
\pi_n(\Map_{\M}(X_{\beta_{i-1}}, \widetilde{X}_{\zeta_i})) &\cong [\Sigma^n X_{\beta_{i-1}}, \widetilde{X}_{\zeta_i}],
\\
\pi_n(\Map_{\M}(X_{\beta_{i}}, \widetilde{X}_{\zeta_i})) &\cong [\Sigma^n X_{\beta_{i}}, \widetilde{X}_{\zeta_i}].
\end{align*}
Thus, the final part of Mayer-Vietoris sequence \ref{MVHspace} of (\ref{plbk}) looks as follows:
$$\xymatrix{\us{i}{\bigoplus}[\Sigma X_{\beta_{i}}, \widetilde{X}_{\beta_i}] \oplus \us{i}{\bigoplus}[\Sigma X_{\zeta_{i}}, \widetilde{X}_{\zeta_i}] \cong \pi_1 \ar[d] \\ \us{i}{\bigoplus}[\Sigma X_{ \beta_{i-1}}, \widetilde{X}_{\zeta_i}] \oplus \us{i}{\bigoplus}[\Sigma X_{\beta_{i}}, \widetilde{X}_{\zeta_i}] \cong \pi_1 \ar[d]^{\partial} \\\;\;\;\;\;\;\;\;\; K \cong \pi_0 \ar[d] \\ \us{i}{\bigoplus}[X_{\beta_{i}}, \widetilde{X}_{\beta_i}] \oplus \us{i}{\bigoplus}[X_{\zeta_{i}}, \widetilde{X}_{\zeta_i}] \cong \pi_0 \ar[d]  \\ \us{i}{\bigoplus}[X_{\beta_{i-1}}, \widetilde{X}_{\zeta_i}] \oplus \us{i}{\bigoplus}[X_{\beta_{i}}, \widetilde{X}_{\zeta_i}] \cong \pi_0, }$$
where $\partial$ is the connecting homomorphism.

Next we would like to have a more explicit description of the terms involved in this exact sequence. Proposition \ref{II}, yields the following isomorphisms:
\begin{align*}
[X_{\beta_i}, \widetilde{X}_{\beta_i}] &\cong \Hom_{\pi_*P}(B^{(i)}(X), B^{(i)}(\widetilde{X})),
\\
[X_{\zeta_i}, \widetilde{X}_{\zeta_i}] &\cong \Hom_{\pi_*P}(Z^{(i)}(X), Z^{(i)}(\widetilde{X})),
\\
[X_{\beta_i}, \widetilde{X}_{\zeta_i}] &\cong \Hom_{\pi_*P}(B^{(i)}(X), Z^{(i)}(\widetilde{X})),
\\
[\Sigma X_{\beta_{i-1}}, \widetilde{X}_{\zeta_i}] &\cong \Hom_{\pi_*P}(B^{(i-1)}(X)[1], Z^{(i)}(\widetilde{X})),
\\
[X_{\beta_{i-1}}, \widetilde{X}_{\zeta_i}] &\cong \Ext^1_{\pi_*P}(B^{(i-1)}(X)[1], Z^{(i)}(\widetilde{X})),\end{align*}
and
\begin{align*}
[\Sigma X_{\beta_{i}}, \widetilde{X}_{\beta_i}] &\cong \Ext^{N-1}_{\pi_*P}(B^{(i)}(X)[N], B^{(i)}(\widetilde{X})),
\\
[\Sigma X_{\zeta_{i}}, \widetilde{X}_{\zeta_i}] &\cong \Ext^{N-1}_{\pi_*P}(Z^{(i)}(X)[N], Z^{(i)}(\widetilde{X})),
\\
[\Sigma X_{\beta_{i}}, \widetilde{X}_{\zeta_i}] &\cong \Ext^{N-1}_{\pi_*P}(B^{(i)}(X)[N], Z^{(i)}(\widetilde{X})).\end{align*}
Since $X \in \K$, the last three Ext-groups vanish and therefore
$$[\Sigma X_{\beta_{i}}, \widetilde{X}_{\beta_i}]=0,\; [\Sigma X_{\zeta_{i}}, \widetilde{X}_{\zeta_i}]=0, \;[\Sigma X_{\beta_{i}}, \widetilde{X}_{\zeta_i}]=0.$$
Consequently, one obtains the desired exact sequence describing $K$:
$$\xymatrix{0 \ar[d] \\ \us{i}{\bigoplus}\Hom_{\pi_*P}(B^{(i-1)}(X)[1], Z^{(i)}(\widetilde{X})) \ar[d]^{\partial} \\ K \ar[d] \\ \us{i}{\bigoplus}\Hom_{\pi_*P}(B^{(i)}(X), B^{(i)}(\widetilde{X})) \oplus \us{i}{\bigoplus} \Hom_{\pi_*P}(Z^{(i)}(X), Z^{(i)}(\widetilde{X})) \ar[d] \\ \us{i}{\bigoplus}\Ext^1_{\pi_*P}(B^{(i-1)}(X)[1], Z^{(i)}(\widetilde{X})) \oplus  \us{i}{\bigoplus}  \Hom_{\pi_*P}(B^{(i)}(X), Z^{(i)}(\widetilde{X})).}$$
As $B^{(i)}(X)$, $B^{(i)}(\widetilde{X})$, $Z^{(i)}(X)$, $Z^{(i)}(\widetilde{X})$ $\in$ $\B[i]$, $i \in \Z/N\Z$, we may rewrite this exact sequence as follows
$$\xymatrix{0 \ar[d] \\ \Hom_{\pi_*P}(B_*(X)[1], Z_*(\widetilde{X})) \ar[d]^{\partial} \\ K \ar[d] \\ \Hom_{\pi_*P}(B_*(X), B_*(\widetilde{X})) \oplus \Hom_{\pi_*P}(Z_*(X), Z_*(\widetilde{X})) \ar[d] \\ \Ext^1_{\pi_*P}(B_*(X)[1], Z_*(\widetilde{X})) \oplus \Hom_{\pi_*P}(B_*(X), Z_*(\widetilde{X})).}$$

Let us now describe
$$L=\Hom_{(\Grmod \pi_*P)^{\underline{1}}}(Z_*(X)  \xto{\iota} C_*(X), Z_*(\widetilde{X})  \xto{\widetilde{\iota}} C_*(\widetilde{X})).$$
First note that any $f \in L$ yields (as we have already mentioned) a diagram
$$\xymatrix{C_*(X) \ar[r]^{\rho} \ar[d]^{f} & B_*(X)[1] \ar[r]^{\lambda[1]} \ar[d]^{\ol{f}[1]} & Z_*(X)[1] \ar[r]^{\iota[1]} \ar[d]^{f[1]} & C_*(X)[1] \ar[d]^{f[1]} \\ C_*(\widetilde{X}) \ar[r]^{\widetilde{\rho}} & B_*(\widetilde{X})[1] \ar[r]^{\widetilde{\lambda}[1]} & Z_*(\widetilde{X})[1] \ar[r]^{\widetilde{\iota}[1]} & C_*(\widetilde{X})[1],}$$
in which the left and the right squares commute. Consider the map $$\sigma \colon L \lra  \Hom_{\pi_*P}(B_*(X), B_*(\widetilde{X})) \oplus \Hom_{\pi_*P}(Z_*(X), Z_*(\widetilde{X})), \;\; \; \;  \sigma(f)=(\ol{f}, f).$$
The kernel of $\sigma$ consists of morphisms of the form $(0,G) \in L$ such that the diagram $$\xymatrix{0 \ar[r] & Z_*(X) \ar[d]^0 \ar[r]^{\iota} & C_*(X) \ar[d]^{G} \ar[r]^{\rho} & B_*(X)[1] \ar[d]^0 \ar[r] & 0 \\ 0 \ar[r] & Z_*(\wt{X}) \ar[r]^{\wt{\iota}} & C_*(\wt{X}) \ar[r]^{\wt{\rho}} & B_*(\wt{X})[1] \ar[r] & 0}$$ commutes, i.e., the kernel can be identified with the subgroup of $\Hom_{\pi_*P}(C_*(X), C_*(\widetilde{X}))$ consisting of those $G \colon C_*(X) \longrightarrow C_*(\wt{X})$ for which $G \iota=0$ and $\wt{\rho}G=0$. It then follows form the snake lemma that $G \in \Ker {\sigma}$ if and only if $G$ can be factored as
$$\xymatrix{C_*(X) \ar[r]^{\rho} & B_*(X)[1] \ar[r]^{g} & Z_*(\wt{X}) \ar[r]^{\wt{\iota}} & C_*(\wt{X})}$$
for some $g$. Since $\rho$ is an epimorphism and $\wt{\iota}$ a monomorphism, $\Ker \sigma$ is
isomorphic to $$\Hom_{\pi_*P}(B_*(X)[1], Z_*(\widetilde{X}))$$ and we get an exact sequence
$$\xymatrix{0 \ar[r] & \Hom_{\pi_*P}(B_*(X)[1], Z_*(\widetilde{X})) \ar[r]^-{b} & L  \ar[r]^-{\sigma} & \Hom_{\pi_*P}(B_*(X), B_*(\widetilde{X})) \oplus \Hom_{\pi_*P}(Z_*(X), Z_*(\widetilde{X})),}$$
where
$$b(g)=(0, \wt{\iota}g\rho).$$
Then, in order to give a more precise description of $L$, we would like to examine $\Imm \sigma$. For this purpose, observe that $$\xymatrix{S \colon \; 0 \ar[r] &  Z_*(X) \ar[r]^{\iota} &  C_*(X) \ar[r]^{\rho} & B_*(X)[1] \ar[r] & 0}$$ and $$\xymatrix{\wt{S} \colon \; 0 \ar[r] & Z_*(\wt{X}) \ar[r]^{\wt{\iota}} &  C_*(\wt{X}) \ar[r]^{\wt{\rho}} & B_*(\wt{X})[1] \ar[r] & 0}$$ represent elements of $\Ext^1_{\pi_*P}(B_*(X)[1], Z_*(X))$ and $\Ext^1_{\pi_*P}(B_*(\wt{X})[1], Z_*(\wt{X}))$, respectively. Define a homomorphism
$$\xymatrix{\Hom_{\pi_*P}(B_*(X), B_*(\widetilde{X})) \oplus \Hom_{\pi_*P}(Z_*(X), Z_*(\widetilde{X})) \ar[r]^-{\nu} & \Ext^1_{\pi_*P}(B_*(X)[1], Z_*(\widetilde{X})),}$$
$$\nu(u,w)=u[1]^*(\wt{S})-w_*(S).$$
We claim that $\Imm \sigma=\Ker \nu$. Indeed, $(u,w) \in \Imm \sigma$ if and only if there exists $h \colon \; C_*(X) \longrightarrow C_*(\wt{X})$ such that the diagram $$\xymatrix{0 \ar[r] & Z_*(X) \ar[d]^w \ar[r]^{\iota} & C_*(X) \ar[d]^{h} \ar[r]^{\rho} & B_*(X)[1] \ar[d]^{u[1]} \ar[r] & 0 \\ 0 \ar[r] & Z_*(\wt{X}) \ar[r]^{\wt{\iota}} & C_*(\wt{X}) \ar[r]^{\wt{\rho}} & B_*(\wt{X})[1] \ar[r] & 0}$$ commutes. But the latter is the case if and only if $w_*(S)=u[1]^*(\wt{S})$, i.e., $(u,w) \in \Ker \nu$. Hence, one concludes that there is an exact sequence
$$\xymatrix{0 \ar[d] \\ \Hom_{\pi_*P}(B_*(X)[1], Z_*(\widetilde{X})) \ar[d]^-{b} \\ L \ar[d]^-{\sigma} \\\Hom_{\pi_*P}(B_*(X), B_*(\widetilde{X})) \oplus \Hom_{\pi_*P}(Z_*(X), Z_*(\widetilde{X})) \ar[d]^-{\nu} \\ \Ext^1_{\pi_*P}(B_*(X)[1], Z_*(\widetilde{X})).}$$

The next step is to compare the exact sequences describing $K$ and $L$. In other words, we will now check that the diagram
\begin{align} \begin{split}\label{dididiag}
\xymatrix{0 \ar[d] & & 0 \ar[d] \\
\Hom_{\pi_*P}(B_*(X)[1], Z_*(\widetilde{X})) \ar[d]^{\partial} \ar@{=}[rr] && \Hom_{\pi_*P}(B_*(X)[1], Z_*(\widetilde{X})) \ar[d]^{b} \\
K \ar[d] \ar[rr]^q && L \ar[d]
\\
\text{$\begin{matrix} \Hom_{\pi_*P}(B_*(X), B_*(\widetilde{X})) \\ \oplus \\ \Hom_{\pi_*P}(Z_*(X), Z_*(\widetilde{X})) \end{matrix}$}
\ar[d] \ar@{=}[rr] &&
\text{$\begin{matrix} \Hom_{\pi_*P}(B_*(X), B_*(\widetilde{X})) \\ \oplus \\ \Hom_{\pi_*P}(Z_*(X), Z_*(\widetilde{X})) \end{matrix}$} \ar[d]  \\
\text{$\begin{matrix} \Ext^1_{\pi_*P}(B_*(X)[1], Z_*(\widetilde{X})) \\ \oplus \\ \Hom_{\pi_*P}(B_*(X), Z_*(\widetilde{X})) \end{matrix}$} \ar[rr]^{\pr} && \Ext^1_{\pi_*P}(B_*(X)[1], Z_*(\widetilde{X})) &}\end{split} \end{align}
commutes. It is clear that the middle square in this diagram is commutative. The commutativity of the lower one is an immediate formal consequence of the construction of the Adams spectral sequence \ref{ASS}. It thus remains to check that
$$q\partial=b$$
which needs a detailed verification.

Fix $j$ $\in$ $\Z/N\Z$ and take $\theta$ $\in$ $\Hom_{\pi_*P}(B^{(j-1)}(X)[1], Z^{(j)}(\widetilde{X}))$. In order to compute $\partial(\theta)$, consider the isomorphisms
$$\pi_1(\Map_{\M}(X_{\beta_{j-1}},\widetilde{X}_{\zeta_j})) \cong [\Sigma X_{\beta_{j-1}}, \widetilde{X}_{\zeta_j}] \cong \Hom_{\pi_*P}(B^{(j-1)}(X)[1], Z^{(j)}(\widetilde{X})).$$
A representative of the homotopy class in $[\Sigma X_{\beta_{j-1}}, \widetilde{X}_{\zeta_j}]$ to which $\theta$ corresponds will be denoted by
$$\theta \colon S^1 \wedge X_{\beta_{j-1}} \longrightarrow \wt{X}_{\zeta_j},$$
abusing notations. Let $\theta_1$ denote the composite
$$\xymatrix{\theta(\pi_+ \wedge 1) \colon I_+  \wedge X_{\beta_{j-1}} \ar[r]^-{\pi_+ \wedge 1} & S^1 \wedge X_{\beta_{j-1}} \ar[r]^-{\theta} & \wt{X}_{\zeta_j},}$$
where $\pi_+ \colon I_+ \xto{} S^1$ is the projection. As noted above,
$$ X_{\beta_{j-1}} \vee X_{\beta_{j}} \xto{(k_j,\; l_j)} X_{\zeta_j}$$
is a cofibration. Besides, the map $\inn_1 \colon S^0 \xto{} I_+$ is an  acyclic cofibration of pointed simplicial sets (see \ref{sset}). Therefore, the pushout product \cite[2.3.11]{GJ99}
$$\inn_1  \square (k_j,\; l_j) \colon (I_+ \wedge (X_{\beta_{j-1}} \vee X_{\beta_{j}})) \vee_{S^0 \wedge (X_{\beta_{j-1}} \vee X_{\beta_{j}})} S^0 \wedge X_{\zeta_j} \longrightarrow I_+ \wedge X_{\zeta_j}$$
is an acyclic cofibration in $\M$. Then, since $\wt{X}_{\zeta_j}$ is fibrant, there exists a morphism $$\Theta \colon I_+ \wedge X_{\zeta_j} \longrightarrow \wt{X}_{\zeta_j},$$ making the diagram
\begin{align} \begin{split} \xymatrix{(I_+ \wedge (X_{\beta_{j-1}} \vee X_{\beta_{j}})) \vee_{S^0 \wedge (X_{\beta_{j-1}} \vee X_{\beta_{j}})} S^0 \wedge X_{\zeta_j}  \ar[rr]^-{(\theta_1,\ast, \ast)} \ar[d]_{\inn_1  \square (k_j,\; l_j)} & & \wt{X}_{\zeta_j} \\ I_+ \wedge X_{\zeta_j} \ar[urr]_\Theta & &}\end{split} \end{align}
commute. This morphism together with the map $\inn_0 \colon S^0 \xto{} I_+$ (see \ref{sset}) gives a composition $\theta_0= \Theta (\inn_0 \wedge 1) \colon X_{\zeta_j} \cong S^0 \wedge X_{\zeta_j} \longrightarrow \wt{X}_{\zeta_j}$. Clearly, one has the equalities $\theta_0k_j=\ast$ and $\theta_0l_j=\ast$. Consequently, there is $\wh{\theta} \in \Hom_{\M^{\C_N}}(X, \wt{X}) \cong (\Map_{\M^{\C_N}}(X, \wt{X}))_0$ having $\theta_0$ as the $\zeta_j$-th component and zero as all the other components. By adjunction, the commutative diagram $(5)$ is equivalent to the commutative square
$$\xymatrix{ S^0 \ar[rr]^-{\ast} \ar[d]_{\inn_1} & & \Map_{\M}(X_{\zeta_j}, \widetilde{X}_{\zeta_j}) \ar[d]^{(k_j^{\ast}, l_j^{\ast})} \\ I_+ \ar[rr]_-{(\theta_1, \ast)} \ar[urr]^-{\Theta} & & \Map_{\M}(X_{\beta_{i-1}}, \widetilde{X}_{\zeta_i}) \times \Map_{\M}(X_{\beta_i}, \widetilde{X}_{\zeta_i}). }$$
Therefore, by construction of $\partial$ (see Lemma \ref{dMV}), we conclude
$$\partial(\theta)=[\wh{\theta}] \in \pi_0(\Map_{\M^{\C_N}}(X,\widetilde{X}))=[X, \wt{X}].$$
Thus one gets the desired computation of $\partial(\theta)$.

Next we show that
$$q\partial(\theta)=b(\theta),$$
i.e.,
$$q([\wh{\theta}])=b(\theta).$$
For this purpose, we explicitly identify both sides of this equation and then compare them. In order to do so, note once again that the cone construction of \ref{cone} is functorial on $(\M_{\cof})^{\underline{1}}$ and hence we have a commutative diagram
$$\xymatrix{ X_{\beta_{j-1}} \ar[d]_{\ast} \ar[r]^{k_j} & X_{\zeta_j} \ar[d]^{\theta_0} \ar[r]^-{\xi} & \Cone(k_j) \ar[r]^{\delta} \ar@{-->}[d]^{\ol{\theta_0}} & S^1 \wedge X_{\beta_{j-1}} \ar[d]^{\ast} \\ \wt{X}_{\beta_{j-1}} \ar[r]^{\wt{k}_j} & \wt{X}_{\zeta_j} \ar[r]^-{\wt{\xi}}  & \Cone(\wt{k}_j) \ar[r]^{\wt{\delta}} & S^1 \wedge \wt{X}_{\beta_{j-1}}.}$$
It follows from the universal property of pushout that there exists $\eta \colon \Cone(k_j) \longrightarrow \wt{X}_{\zeta_j}$ making the diagram
$$\xymatrix{ X_{\zeta_j} \ar[r]^-{\xi} \ar[dr]_{\theta_0} & \Cone(k_j) \ar[d]^{\eta} & (I,0) \wedge X_{\beta_{j-1}} \ar[l]_-{\mu} \ar[dl]^{\ast} \\ &  \wt{X}_{\zeta_j} & }$$
commute (since $\theta_0k_j=\ast$) and that ${\ol{\theta_0}}=\wt{\xi} \eta$. Then ($\theta_0$ is null-homotopic via $\Theta$)
$$q([\wh{\theta}])=(0, \pi_*\ol{\theta_0})=(0,  \pi_*(\wt{\xi}) \pi_*(\eta))$$
and
$$b(\theta)=(0, \pi_*(\wt{\xi})\pi_*(\theta) \pi_*(\delta))$$
by definitions of $q$ and $b$, respectively. Thus, to verify $q([\wh{\theta}])=b(\theta)$, it suffices to prove that $\eta$ is homotopic to $\theta \delta$. This is done explicitly. Indeed, recall that the simplicial set $I \times I$ is generated by the $2$-simplices $h_0=((0,0,1),(0,1,1))$ and $h_1=((0,1,1),(0,0,1))$ satisfying the relation $\partial_1h_0=\partial_1h_1$. Therefore, there is $\wh{F} \colon I \times I \longrightarrow \Map_{\M}(X_{\beta_{j-1}}, \widetilde{X}_{\zeta_j})$ with $\wh{F}(h_0)= \wh{F}(h_1)=s_0\theta_1$ ($\theta_1$ is considered as an element of ($\Map_{\M}(X_{\beta_{j-1}}, \widetilde{X}_{\zeta_j}))_1$ (see \ref{same letter})). It is easily seen that $\wh{F}\vert_{I \times 0}=\ast$, $\wh{F}\vert_{0 \times I}=\ast$, $\wh{F}\vert_{I \times 1}=\theta_1$ and $\wh{F}\vert_{1 \times I}=\theta_1$ (see \ref{sset}). Hence $\wh{F}$ induces a simplicial map
$$\xymatrix{I_+ \wedge (I,0) \ar[r]^-{\wt{F}} & \Map_{\M}(X_{\beta_{j-1}}, \widetilde{X}_{\zeta_j})}$$
with $\wt{F}\vert_{0 \wedge (I,0)}=\ast$, $\wt{F}\vert_{1 \wedge (I,0)}=\theta_1$, and $\wt{F}\vert_{I_+ \wedge 1}=\theta_1$ (see \ref{sset}). Let
$$\xymatrix{I_+ \wedge ((I,0) \wedge X_{\beta_{j-1}}) \ar[r]^-{F} & \widetilde{X}_{\zeta_j}}$$
be the adjoint of $\wt{F}$. Clearly, $ F \vert_{0 \wedge (I,0)  \wedge X_{\beta_{j-1}} }=\ast$, $ F \vert_{ 1 \wedge (I,0) \wedge X_{\beta_{j-1}} }=\theta(\pi \wedge 1 )= \theta \delta \mu$ (where $\pi \colon I \xto{} S^1$ is the projection), and  $F \vert_{I_+ \wedge 1 \wedge X_{\beta_{j-1}}}=\theta_1$. The latter shows that the diagram
$$\xymatrix{I_+ \wedge X_{\beta_{j-1}} \ar[r]^{1 \wedge k_j} \ar[d] & I_+ \wedge X_{\zeta_j} \ar[d]^{\Theta} \\ I_+ \wedge ((I,0) \wedge X_{\beta_{j-1}}) \ar[r]^-{F} & \widetilde{X}_{\zeta_j}}$$
commutes. (Note that the left vertical arrow in this diagram comes from the morphism
$$\xymatrix{X_{\beta_{j-1}} \cong S^0 \wedge X_{\beta_{j-1}} \ar[r]^-{\incl \wedge 1} & (I,0) \wedge X_{\beta_{j-1}}.})$$
On the other hand, by \ref{cone}, we have a pushout square
$$\xymatrix{I_+ \wedge X_{\beta_{j-1}} \ar[r]^{1 \wedge k_j} \ar[d] & I_+ \wedge X_{\zeta_j} \ar[d]^{1 \wedge \xi} \\ I_+ \wedge ((I,0) \wedge X_{\beta_{j-1}})\ar[r]^-{1 \wedge \mu} & I_+ \wedge \Cone(k_j)}$$ since $I_+ \wedge -$ is left adjoint.
Consequently, there exists $H \colon I_+ \wedge \Cone(k_j) \longrightarrow \widetilde{X}_{\zeta_j}$ such that the diagram
$$\xymatrix{ I_+ \wedge X_{\zeta_j} \ar[r]^-{1 \wedge \xi} \ar[dr]_{\Theta} & I_+ \wedge \Cone(k_j) \ar[d]^-{H} & I_+ \wedge ((I,0) \wedge X_{\beta_{j-1}}) \ar[l]_-{1 \wedge \mu} \ar[dl]^{F} \\ &  \wt{X}_{\zeta_j} &}$$ commutes. A trivial computation shows that H is a left homotopy from $\eta$ to $\theta \delta$. Hence $q\partial=b$.

Thus (\ref{dididiag}) commutes. Now a diagram chase shows that (\ref{eseq}) is exact, and this completes the proof of \ref{V}. \end{prf}

Next we examine the image of $\Q\vert_{ \K} \colon \K \longrightarrow \Mod \pi_*P$. By the construction of $\Q$ and the definition of $\K$, we have
$$\proj (Z_*(X)=\Ker d) < N-1,\;\;\;\;\;\;\;\;\; \proj (B_*(X)= \Img d) < N-1,\;X \in \K.$$
These key properties completely determine the essential image of $\Q\vert_{ \K}$. Indeed, one has

\begin{prop}\label{VI} A differential graded $\pi_*P$-module $(C_*, d)$ is in the essential image of $\Q\vert_{ \K}$ if and only if $$\proj  Z_*<N-1\;\;\;\;\;\;\; \text{and} \;\;\;\;\;\;\;  \proj  B_*<N-1,$$
where $Z_* = \Ker d$ and $B_* = \Imm d.$
\end{prop}

\begin{prf} We have just mentioned that the ``only if'' part of this proposition is valid. Let us check that the ``if'' part holds as well.

Since $\pi_*P$ is $N$-sparse, we have the decompositions (see \ref{spectral}) $$Z_*= Z^{(0)} \oplus Z^{(1)} \oplus \ldots \oplus Z^{(N-1)},$$ $$B_*= B^{(0)} \oplus B^{(1)} \oplus \ldots \oplus B^{(N-1)},$$ $$C_*= C^{(0)} \oplus C^{(1)} \oplus \ldots \oplus C^{(N-1)},$$ where $C^{(i)}, B^{(i)}, Z^{(i)} \in \B[i]$, $i \in \Z/N\Z$. By \ref{II} (i) there exist $X_{\zeta_i}, X_{\beta_i} \in \E_i$ such that
$$Z^{(i)} \cong \pi_*X_{\zeta_i} \;\;\;\;\;\;\; \text{and} \;\;\;\;\;\;\; B^{(i)} \cong \pi_*X_{\beta_i}$$
for all $i \in \Z/N\Z$. We fix these isomorphisms once and for all. Without loss of generality one may assume that $X_{\beta_i}$ is cofibrant and $X_{\zeta_i}$ is bifibrant in $\M$, $i \in \Z/n\Z$. It then follows from Proposition \ref{II} (i) that there exist morphisms $l_i \colon X_{\beta_i} \longrightarrow X_{\zeta_i}$ in $\M$ such that the diagrams
$$\xymatrix{ \pi_*X_{\beta_i} \ar[r]^{\pi_*l_i} \ar[d]^{\cong} & \pi_*X_{\zeta_i} \ar[d]^{\cong} \\ B^{(i)} \; \ar@{>->}[r] & Z^{(i)}}$$
commute. Further, for any $i \in \Z/N\Z$, consider the extension $$\xymatrix{ E_i \colon 0 \ar[r] & Z^{(i)} \ar[r] & C^{(i)} \ar[r] & B^{(i-1)}[1] \ar[r] & 0.}$$ Under Isomorphism \ref{II}(ii), the extension class of $E_i$ corresponds to some element in $[X_{ \beta_{i-1}}, X_{\zeta_i}]$ which is the homotopy class of some morphism $k_i \colon X_{\beta_{i-1}} \longrightarrow  X_{\zeta_i}$ in $\M$. Consequently, we get a diagram
$$\xymatrix{X \colon & X_{\zeta_0}
    \ar@{<-}[drrrrrrrr]%
    |<<<<<<<<<<<<<<<<<<{\text{\Large \textcolor{white}{$\blacksquare$}}}%
    |<<<<<<<<<<<<<<<<<<<<<<<{\text{\Large \textcolor{white}{$\blacksquare$}}}
    |<<<<<<<<<<<<<<<<<<<<<<<<<<<<<<<<<<<<<<<{\text{\Large \textcolor{white}{$\blacksquare$}}}
    |>>>>>>>>>>>>>>>>>>>>>>>>>>>>>>>>>>>>>>>>{\text{\Large \textcolor{white}{$\blacksquare$}}}
    |>>>>>>>>>>>>>>>>>>>>>>>>>{\text{\Large \textcolor{white}{$\blacksquare$}}}
    |>>>>>>>>>>>>>>>>>>{\text{\Large \textcolor{white}{$\blacksquare$}}}
& & X_{\zeta_1} & & \ldots & & X_{\zeta_{N-2}} & & X_{\zeta_{N-1}} \\  & X_{\beta_0} \ar[u]^{l_0} \ar[urr]^{k_1} & & X_{\beta_1} \ar[u]_{l_1} \ar[urr] & & \ldots \ar[urr] & & X_{\beta_{N-2}} \ar[u]^{l_{N-2}} \ar[urr]^{k_{N-1}} & & X_{\beta_{N-1}}, \ar[u]_{l_{N-1}} }$$
i.e., an object in $\M^{\C_N}$. Without loss of generality we may assume that $X \in (\M^{\C_N})_{\cof}$ (otherwise we could replace $X$ cofibrantly). Besides, $\proj  Z_*$ and $\proj  B_*$ are less than $N-1$. Therefore, $X$ is an object of $\K$.  Finally, the construction of $\Q$ and Remark \ref{III} immediately imply that $\Q(X)$ is isomorphic $C_*$. \end{prf}

In view of \ref{V}, one gets

\begin{coro} \label{Qequicase} Let $\M$ be a stable model category and $P$ a compact generator of $\Ho(\M)$. Suppose that the graded ring $\pi_*P$ is $N$-sparse and $\gl \pi_*P < N-1$. Then the functor
$$\Q \colon \L \longrightarrow \Mod \pi_*P $$
is an equivalence of categories.
\end{coro}

We conclude this subsection with

\begin{prop}\label{VII} For any $X \in \L$ there is a natural isomorphism of graded $\pi_*P$-modules
$$H_*(\Q(X)) \cong \pi_*(\Hocolim X).$$
\end{prop}

\begin{prf}
As noted in Remark \ref{IV}, $\L $ is the localization of $\L' \subseteq (\M^{\C_N})_{\cof}$ at the class of weak equivalences. Let
$$\Q'=\Q \gamma \colon \L' \longrightarrow \Mod \pi_*P,$$
where $\gamma \colon \L' \longrightarrow \L$  is the localization functor (which is the identity on objects!). Clearly, to prove the proposition, it suffices to construct a natural isomorphism
$$H_*(\Q'(X)) \cong \pi_*(\colim X), \;\;\;\;\;\;\;\;\; X \in \L'.$$
Recall that $X_{\beta_{i}}$ and $X_{\zeta_i}$ are cofibrant and
$$ X_{\beta_{i-1}} \vee X_{\beta_{i}} \xto{(k_i,\; l_i)} X_{\zeta_i}$$
is a cofibration for all $ i \in \Z/N\Z$ (since $X$ is cofibrant in $\M^{\C_N}$ (see \ref{exmpV})). This implies that the coequalizer diagram
\begin{align}\begin{split}\label{hocoeq}\xymatrix{\vee_{i}X_{\beta_{i}} \ar@<0.5ex>[r]^k \ar@<-0.5ex>[r]_l & \vee_{i}X_{\zeta_i} \ar[r] & \colim X,}\end{split}\end{align}
where $k=\vee_{i}k_i$ and $l=\vee_{i}l_i$, is in fact a homotopy coequalizer. Hence, we get a long exact sequence (see \ref{htpycoeq})
\begin{multline*} \ldots \xto{}  \pi_*(\colim X)[-1] \xto{\partial}  \oplus_{i} \pi_*X_{\beta_{i}} \xto{\pi_*l-\pi_*k} \oplus_{i} \pi_*X_{\zeta_{i}} \xto{} \pi_*(\colim X) \xto{\partial}  \oplus_{i} \pi_*X_{\beta_{i}}[1] \xto{} \ldots .\end{multline*}
Since $X \in \L'$, $\pi_*l=\oplus_{i}\pi_*l_i$ is a monomorphism and $\pi_*k=\oplus_{i}\pi_*k_i$ is zero. Therefore, the connecting homomorphisms of this long exact sequence vanish and we get a short exact sequence
$$ \xymatrix{0 \ar[r] & \oplus_i \pi_*X_{\beta_{i}} \ar[rr]^{\oplus_{i}\pi_*l_i} & & \oplus_i \pi_*X_{\zeta_{i}} \ar[r] & \pi_*(\colim X) \ar[r] & 0,}$$
i.e.,
$$ \xymatrix{0 \ar[r] & B_*(X) \ar[r]^{\lambda} & Z_*(X) \ar[r] & \pi_*(\colim X) \ar[r] & 0.}$$
This short exact sequence shows that there is a natural isomorphism
$$ Z_*(X) / B_*(X) \cong \pi_*(\colim X).$$
The short exact sequence (\ref{se}) together with the definition of the differential of $C_*(X)$ gives
$$H_*(\Q'(X))=H_*(C_*(X)) \cong Z_*(X)/B_*(X).$$ Thus, we get a natural isomorphism
$$H_*(\Q'(X)) \cong \pi_*(\colim X).$$ \end{prf}

\vspace{0.75em}

\subsection{Proof of \ref{I} (i)}~\label{R}

\vspace{0.75em}

\noindent Our aim is to construct a functor
$$\R \colon \D(\pi_*P) \longrightarrow \Ho(\M),$$
and a natural isomorphism $\pi_* \circ \R \cong H_*$.

It follows from Proposition \ref{VI} that the essential image $\Q(\K)$ of $\Q\vert_{ \K}$ contains all differential graded $\pi_*P$-modules which have underlying projective graded $\pi_*P$-modules. By \cite[2.2.5]{H97},  it in particular contains all cofibrant DG $\pi_*P$-modules. Therefore, the derived category $\D(\pi_*P)$ is equivalent to the localization $\Q(\K)[\V^{-1}]$, where $\V$ is the class of quasi-isomorphisms in $\Q(\K)$. Thus, in order to construct $\R$, we need a functor
$$\R' \colon \Q(\K) \longrightarrow \Ho(\M)$$
which sends quasi isomorphisms of DG modules to isomorphisms. By Proposition \ref{V}, $$\Q\vert_{ \K} \colon \K \longrightarrow \Q(\K)$$ is an equivalence of categories and so we can choose an inverse
$(\Q\vert_{ \K})^{-1} \colon \Q(\K) \longrightarrow \K.$ Define $\R'$ to be the following composite
$$\xymatrix{ \Q(\K) \ar[rr]^-{(\Q\vert_{ \K})^{-1}} & & \K \subseteq \Ho(\M^{\C_N}) \ar[rr]^-{\Hocolim} & & \Ho(\M).}$$
Then, Proposition \ref{VII} immediately yields a natural isomorphism
$$\pi_*(\R'(M)) \cong H_*M, \;\;\;\;\;\;\;\;\; M \in \Q(\K).$$ In other words, the diagram
\begin{align}\begin{split}\label{diag1}\xymatrix{\Q(\K) \ar[rr]^{H_*} \ar[dr]_{\R'} & & \Grmod{\pi_*P} \\ & \Ho(\M) \ar[ur]_{\pi_*} &}\end{split}\end{align}
commutes up to a natural isomorphism. This shows that $\R'$ sends quasi-isomorphisms to isomorphisms (since $\pi_*$ reflects isomorphisms).  Consequently, $\R'$ factors through the localization, i.e., there exists $$\R \colon \D(\pi_*P) \longrightarrow \Ho(\M)$$ such that the diagram
\begin{align}\begin{split}\label{diag2}\xymatrix{\Q(\K) \ar[rr]^{\R'} \ar[dr]_{\Lambda} & & \Ho(\M) \\ & \D(\pi_*P) \ar[ur]_{\R}, &}\end{split}\end{align}
where $\Lambda$ is the localization functor, commutes up to a natural isomorphism. The universal property of localization and diagrams (\ref{diag1}) and (\ref{diag2}) and give the desired result.\;\;\;\;\;\;\;\;\;\;\;\;\;\;\;\;\;\;\;\;\;\;\;\;\;\;\;\;\;\;\;$\Box$

\subsection{Preservation of suspensions}

\begin{prop}\label{commshift} Let $\M$ be a stable model category and $P$ a fixed compact generator of $\Ho(\M)$. Suppose that the graded ring $\pi_*P$ is $N$-sparse and the (right) graded global homological dimension of $\pi_*P$ is less than $N$. Then the functor
$$\R \colon \D(\pi_*P) \longrightarrow \Ho(\M)$$
commutes with suspensions. \end{prop}

\begin{prf} We consider the functor
$$ (-)^\# \colon \L \longrightarrow \L$$
sending $X \in \L$ to $X^\#$ which is defined by
$$ X^\#_{ \beta_i} = S^1 \wedge X_{\beta_{i-1}}, \;\;\;\;\;\; X^\#_{ \zeta_i} = S^1 \wedge X_{\zeta_{i-1}},$$
$$ k^\#_i= 1 \wedge k_{i-1}, \;\;\; l^\#_i= 1 \wedge l_{i-1}, \;\;\; i \in \Z/N\Z.$$
There is a natural isomorphism
$$ \Hocolim(X^\#) \cong \Sigma \Hocolim X, \;\; X \in \L$$
in $\Ho(\M)$. On the other hand, by construction of $\Q$, we have a natural isomorphism
$$ \Q(X^\#) \cong \Q(X)[1], \;\; X \in \L$$
(see \ref{wanacvleba}). Consequently, for any $M \in \Q(\K)$, there is a natural isomorphism
$$ \Q^{-1}(M[1]) \cong (\Q^{-1}(M))^\#.$$
By construction of $\R'$, one gets
\begin{align*}\begin{split} \R'(M[1]) = \Hocolim (\Q^{-1}(M[1])) \cong \Hocolim((\Q^{-1}(M))^\#) \cong \\ \Sigma \Hocolim (\Q^{-1}(M))= \Sigma \R'(M),\end{split} \end{align*}
whence we have a natural isomorphism
$$\R(M[1]) \cong \Sigma \R(M), \;\; M \in \D(\pi_*P).$$ \end{prf}

\section{Further properties of $\R$ and examples}\label{appl}

In this section we prove that the functor $\R$ need not be a derived functor of a Quillen functor. Besides, we give non-trivial examples where \ref{I} can be applied.

First let us briefly recall

\vspace{0.75em}

\subsection{Derived mapping spaces}~\label{csframes}

\vspace{0.75em}

\noindent For a pointed simplicial model category $\M$, one has the mapping space functor
$$\Map_{\M}(-,-) \colon \M^{op} \times \M \longrightarrow \sSet_*.$$
The simplicial model category axioms imply that this functor is in fact a Quillen bifunctor (see \cite[4.2.1]{H99}). Therefore, we get the derived mapping space functor
$$\Map_{\Ho(\M)}(-,-) \colon \Ho(\M)^{op} \times \Ho(\M) \longrightarrow \Ho(\sSet_*).$$
If a model category $\M$ is not simplicial, then one does not have mapping spaces on the model level. However, one is still able to define derived mapping spaces. This is done in \cite[5.4]{H99} as follows.

Let $\M$ be a model category and $\M^\Delta$ denote the category of cosimplicial objects in $\M$. Consider the functor
$$\Ev_0 \colon \M^\Delta \longrightarrow \M$$
which sends $X^{\bullet} \in \M^\Delta$ to the $0$-th space $X^{\bullet}[0]$. One can easily see that $\Ev_0$ has both a right and left adjoint. The right adjoint
$${\rr}^{\bullet} \colon \M \longrightarrow \M^\Delta$$
is the constant cosimplicial object functor. The left adjoint
$${\cop}^{\bullet} \colon \M \longrightarrow \M^\Delta$$
sends $X \in \M$ to the cosimplicial object whose $n$-th space is the $n+1$-fold coproduct of $X$ and whose structure morphisms are obtained from the coproduct inclusions and fold maps. The functors ${\cop}^{\bullet}$ and ${\rr}^{\bullet}$ come with a canonical natural transformation
$$\cop^\bullet \longrightarrow \rr^\bullet$$
which is the identity in degree $0$ and the fold map in higher degrees.

Note that there is a model structure on $\M^\Delta$, called the Reedy model structure (see \cite[5.2]{H99} for details), and that the adjunctions $({\cop}^{\bullet}, \Ev_0)$, $(\Ev_0, {\rr}^{\bullet})$ are in fact Quillen adjunctions.

\begin{defi}[{\cite[5.2.7]{H99}}]\label{frames}  Suppose $\M$ is a model category and $X$ an object of $\M$. A cosimplicial frame on $X$ is a cosimplicial object $X^{\bullet} \in \M^\Delta$ together with a factorization
$$\xymatrix{\cop^{\bullet} X \ar[r] & X^\bullet \ar[r] & \rr^\bullet X}$$
of the canonical map $$\cop^\bullet X \longrightarrow \rr^\bullet X$$ into a cofibration in $\M^\Delta$ followed by a weak equivalence (which is an isomorphism in degree zero). \end{defi}

The existence of such frames is shown in \cite[5.2.8]{H99}.

\;

\begin{defi}\label{dhoms} Let $X,Y \in \Ho(\M)$. The derived mapping space $\Map_{\Ho(\M)}(X,Y)$ is defined by
$$\Map_{\Ho(\M)}(X,Y) = \Hom_{\M}( X^{\bullet}, Y^f) \in \sSet,$$
where $X^{\bullet}$ is a cosimplicial frame of $X$ and $Y^f$ a fibrant replacement of $Y$. \end{defi}

It is proved in \cite[5.4]{H99} that the derived mapping space is a well-defined object of $\Ho(\sSet)$, i.e., it does not depend, up to homotopy, on the choice of the cosimplicial frame $X^\bullet$ and the fibrant replacement $Y^f$. Moreover, we have a functor
$$\Map_{\Ho(\M)}(-,-) \colon \Ho(\M)^{op} \times \Ho(\M) \longrightarrow \Ho(\sSet).$$
If $\M$ is pointed, then $\Map_{\Ho(\M)}(X,Y)$ comes with a natural base point and it is a well-defined object of $\Ho(\sSet_*)$.

\begin{prop}[{\cite[5.6.2]{H99}}]\label{compare} Suppose
$$\xymatrix{ F \colon \M \ar@<0.5ex>[r] & \N \cocolon G \ar@<0.5ex>[l]}$$
is a Quillen adjunction, where $F$ is a left adjoint. Then, for any $X \in \Ho(\M)$ and $Y \in \Ho(\N)$, there is a natural weak equivalence of derived mapping spaces
$$ \Map_{\Ho(\N)}(\textbf{L} FX,Y) \simeq \Map_{\Ho(\M)}(X,\textbf{R}GY).$$ \end{prop}

\begin{remk}\label{dgmaps} Let $A$ be a differential graded ring. The category $\Mod A$ with the model structure \ref{exmpIII} is not a simplicial model category. However, one can apply \ref{dhoms} to $\Mod A$ and get derived mapping spaces. Note that since $\Mod A$ is an additive category, $\Map_{\D(A)}(X,Y)$ is in fact a simplicial abelian group for any $X,Y \in \D(A)$. It is well known that every simplicial abelian group has the homotopy type of a product of Eilenberg-MacLane spaces (see e.g. \cite[III.2.18]{GJ99}). In particular, $\Map_{\D(A)}(X,Y)$ is homotopy equivalent to a product of Eilenberg-MacLane spaces for any $X,Y \in \D(A)$. \end{remk}

\begin{remk}\label{prfframings} There is an analog of the pushout product axiom for cosimplicial frames (see \cite[5.4.1 and 5.4.2]{H99}). Using this and the derived mapping spaces it is not difficult to see that the proofs given in Sections \ref{yv. mtavari}, \ref{advgamoy} and \ref{intgamoy} can be applied to stable model categories without a simplicial enrichment. \end{remk}

\vspace{0.75em}

\subsection{The cases when $\R$ is not a derived functor}~ \label{Rnonderived}

\vspace{0.75em}

\noindent The following shows that \ref{compare} is rather useful.

\begin{theo}\label{noderived} Let $\M$ be a stable model category and $P$ a compact generator of $\Ho(\M)$. Suppose that the graded ring $\pi_*P$ is $N$-sparse and $\gl \pi_*P < N$ (see \ref{grdims}). Suppose further that the derived mapping space $\Map_{\Ho(\M)}(P,P)$ is not weakly equivalent to a product of Eilenberg-MacLane spaces. Then the functor
$$\R \colon \D(\pi_*P) \longrightarrow \Ho(\M)$$
is not a derived functor of a Quillen functor. \end{theo}

\begin{prf} Consider the object
$$\xymatrix{Y \colon & P
    \ar@{<-}[drrrr]%
    |<<<<<<<<<<<{\text{\Large \textcolor{white}{$\blacksquare$}}}%
    |<<<<<<<<<<<<<<<<<<<{\text{\Large \textcolor{white}{$\blacksquare$}}}
    |>>>>>>>>>>>>>>>>>>>{\text{\Large \textcolor{white}{$\blacksquare$}}}
    |>>>>>>>>>>>{\text{\Large \textcolor{white}{$\blacksquare$}}}
    |>>>>>>>>>{\text{\Large \textcolor{white}{$\blacksquare$}}}
& \ast & \ldots & \ast & \ast \\ & \ast \ar[u] \ar[ur] & \ast \ar[u] \ar[ur] & \ldots \ar[ur] & \ast \ar[u] \ar[ur] & \ast \ar[u]}$$
in $\M^{\C_N}$. Clearly, $Y \in \K$ (see \ref{V}). By construction of $\Q$,
$$\Q(Y) \cong \pi_*P,$$
where $\pi_*P$ is regarded as a DG $\pi_*P$-module with zero differential. On the other hand, $$\Hocolim Y=\colim Y \cong P.$$ Consequently,
$$\R(\pi_*P) \cong P$$
in $\Ho(\M)$.
Now assume that $\R$ is the derived functor of a left Quillen functor. Then $\R$ has a right adjoint
$$\G \colon \Ho(\M) \longrightarrow \D(\pi_*P)$$
and, in view of \ref{compare}, we get a weak equivalence of derived mapping spaces
$$ \Map_{\Ho(\M)}(\R(\pi_*P),P) \simeq \Map_{\D(\pi_*P)}(\pi_*P,\G(P)),$$
i.e.,
$$ \Map_{\Ho(\M)}(P,P) \simeq \Map_{\D(\pi_*P)}(\pi_*P,\G(P)).$$
By \ref{dgmaps}, the mapping space on the right is a product of Eilenberg-MacLane spaces, a contradiction.

One similarly shows that $\R$ cannot be the derived functor of a right Quillen functor. \end{prf}

\begin{coro}\label{nonderived2} Let $S$ be a symmetric ring spectrum such that the graded homotopy ring $\pi_*S$ is $N$-sparse and $\gl \pi_*S < N$, and the infinite loop space $\Omega^{\infty}S$ is not weakly equivalent to a product of Eilenberg-MacLane spaces. Then the functor
$$\R \colon \D(\pi_*S) \longrightarrow \D(S)$$
is not derived from a Quillen functor.\end{coro}

Further, in view of \ref{compare} and \ref{dgmaps}, one gets the following

\begin{prop} \label{noneq1}Let $\M$ be a stable model category and $P$ a compact generator of $\Ho(\M)$. Suppose that the derived mapping space $\Map_{\Ho(\M)}(P,P)$ is not weakly equivalent to a product of Eilenberg-MacLane spaces. Then there does not exist a zig-zag of Quillen equivalences between the model categories $\M$ and $\Mod \pi_*P$. In particular, if the graded ring $\pi_*P$ is $N$-sparse, $\gl \pi_*P < N$, and the functor
$$\R \colon \D(\pi_*P) \longrightarrow \Ho(\M)$$
is an equivalence of categories (see \ref{shecdoma}), then $\R$ can not be derived from a zig-zag of Quillen equivalences. \end{prop}

\begin{coro}\label{noneq2} Let $S$ be a symmetric ring spectrum. Suppose that the infinite loop space $\Omega^{\infty}S$ is not weakly equivalent to a product of Eilenberg-MacLane spaces. Then there does not exist a zig-zag of Quillen equivalences between the model categories $\Mod S$ and $\Mod \pi_*S$. In particular, if the graded homotopy ring $\pi_*S$ is $N$-sparse, $\gl \pi_*S < N$, and the functor
$$\R \colon \D(\pi_*S) \longrightarrow \D(S)$$
is an equivalence of categories, then $\R$ can not be not derived from a zig-zag of Quillen equivalences.\end{coro}

\subsection{Examples}

\begin{exmp}[Truncated Brown-Peterson spectra]\label{BP} Let $p$ be a prime, $n$ a natural number, and suppose that
$$ n+1 < 2(p-1).$$
Then the truncated Brown-Peterson spectrum $BP \langle n \rangle$ satisfies the hypotheses of \ref{Rcoro}. Indeed, the homotopy ring of $BP \langle n \rangle$ looks as follows
$$\pi_*BP \langle n \rangle \cong \Z_{(p)}[v_1,\ldots,v_n], \;\;\;\;\;\;\;\;\;\;\;\; |v_i|=2(p^i-1).$$
Clearly, $\gl \pi_*BP \langle n \rangle= n+1$ and $\pi_*BP \langle n \rangle$ is $2(p-1)$-sparse.

Note also that $\Omega^{\infty}BP \langle n \rangle$ is not weakly equivalent to a product of Eilenberg-MacLane spaces. Therefore, in view of
\ref{nonderived2}, the functor
$$\R \colon \D(\Z_{(p)}[v_1,\ldots,v_n]) \longrightarrow \D(BP \langle n \rangle)$$
is not a derived functor of a Quillen functor. \end{exmp}

\begin{exmp}[Johnson-Wilson spectra]\label{BP} The Johnson-Wilson spectrum $E(n)$ is obtained from $BP \langle n \rangle$ by inverting the generator $v_n$. In particular, for the homotopy ring of $E(n)$ one has
$$\pi_*E(n) \cong \Z_{(p)}[v_1,\ldots,v_n, {v_n}^{-1}], \;\;\;\;\;\;\;\;\;\;\;\; |v_i|=2(p^i-1).$$
Note that $\gl \pi_*E(n) = n$ and $\pi_*E(n)$ is $2(p-1)$-sparse. Therefore, if
$$ n < 2(p-1),$$
then \ref{Rcoro} can be applied to $E(n)$. Furthermore, since $\Omega^{\infty}E(n)$ is not weakly equivalent to a product of Eilenberg-MacLane spaces, the functor
$$\R \colon \D(\Z_{(p)}[v_1,\ldots,v_n,{v_n}^{-1} ]) \longrightarrow \D(E(n))$$
is not a derived functor of a Quillen functor. \end{exmp}

\section{The one-dimensional case}\label{advgamoy}

In this section we prove \ref{I} (ii).

\subsection{Technical lemmas}

The following well-known proposition immediately follows from Adams spectral sequence \ref{ASS}.

\begin{prop}[Universal coefficient theorem]\label{uct} Suppose $\T$ is a triangulated category with infinite coproducts, $P \in \T$ a compact generator, and let $\pi_* = \Hom_{\T}(P,-)_*$. Then for any $X$ with $\proj  \pi_*X \leq 1$ (in $\Grmod \pi_*P$) and $Y \in \T$, there is a natural short exact sequence
$$ \xymatrix{ 0 \ar[r] & \Ext^{1}_{\pi_*P}( \pi_*X[1], \pi_*Y) \ar[r] & \Hom_{\T}(X,Y) \ar[r]^-{\pi_*} & \Hom_{\pi_*P} (\pi_*X, \pi_*Y) \ar[r] & 0.}$$
In particular, $X$ is isomorphic to $Y$ in $\T$ if and only if $\pi_*X$ and $\pi_*Y$ are isomorphic as graded $\pi_*P$-modules.
\end{prop}

\begin{coro}\label{dashla} Let $\T$, $P$ and $\pi_*$ be as in \ref{uct}, and suppose that $\pi_*P$ is $N$-sparse. Then any $X \in \T$ with $\proj  \pi_*X \leq 1$ splits as follows
$$X \cong \oplus_{i \in \Z/N\Z} X^{(i)},\;\;\;\; \pi_*X^{(i)} \in \B[i],$$
where $\B$ is as in \ref{spectral}. \end{coro}

\begin{defi} Let $\T$ be a triangulated category. One says that a triangle
$$\xymatrix{X \ar[r]^f & Y \ar[r]^g & Z \ar[r]^h & \Sigma X}$$
in $\T$ is antidistinguished if the triangle
$$\xymatrix{X \ar[r]^f & Y \ar[r]^g & Z \ar[r]^{-h} & \Sigma X}$$
is distinguished.

\end{defi}

\begin{lem} \label{trsh} Let $\M$ be a (simplicial) stable model category, $G_0$, $G_1$ cofibrant objects of $\M$, $\omega \colon G_0 \longrightarrow G_1$ a cofibration in $\M$, and $P$ a compact generator of $\Ho(\M)$. Suppose that the graded ring $\pi_*P$ is $N$-sparse and $\gl \pi_*P < N$. Further, assume that $\pi_*G_0, \pi_*G_1$ are in $\B[i]$ for fixed $i \in \Z/N\Z$, $\pi_*G_0$, $\pi_*G_1$ are projective graded $\pi_*P$-modules and $\pi_*\omega$ is a monomorphism. Then the functor
$$\R \colon \D(\pi_*P) \longrightarrow \Ho(\M)$$
sends the antidistinguished triangle (see Remark \ref{algcone})
\begin{align}\begin{split}\label{dim1tr}\xymatrix{\pi_*G_0 \ar[r]^{\pi_*\omega} & \pi_*G_1 \ar[r] & C(\pi_*\omega) \ar[r]^{-\partial} & \pi_*G_0[1]}\end{split}\end{align}
to a triangle isomorphic to the elementary distinguished triangle
$$\xymatrix{G_0 \ar[r]^{\omega} & G_1 \ar[r] & \Cone(\omega) \ar[r] & S^1 \wedge G_0}$$
of $\omega$ (see \ref{cone}).
\end{lem}

\begin{prf} By \ref{VI}, since $\pi_*G_0$, $\pi_*G_1$ are projective, the triangle (\ref{dim1tr}) is
contained in $\Q(\K)$. Let us construct a sequence
\begin{align}\begin{split}\label{dim1seq}\xymatrix{A_0 \ar[r]^{\alpha_0} & A_1 \ar[r]^{\alpha_1} & A_2 \ar[r]^{\alpha_2} & A_3}\end{split}\end{align}
in $\K$ which is mapped to the triangle (\ref{dim1tr}) after applying $\Q$. Define
$$\xymatrix{A_0 \colon & \ast
    \ar@{<-}[drrrrrrr]%
    |<<<<<<<<<<<{\text{\Large \textcolor{white}{$\blacksquare$}}}%
    |<<<<<<<<<<<<<<<<<<<<<{\text{ \Large \textcolor{white}{$\blacksquare$}}}
    |<<<<<<<<<<<<<<<<<<<<<<<<<<<<<<<<<<<<<{\text{\Large \textcolor{white}{$\blacksquare$}}}
    |<<<<<<<<<<<<<<<<<<<<<<<<<<<<<<<<<<<<<<<<<<<<{\text{\Large \textcolor{white}{$\blacksquare$}}}
    |>>>>>>>>>>>>>>>>>>>>>>>>>>>>>>>>>>>>>>{\text{\Large \textcolor{white}{$\blacksquare$}}}
    |>>>>>>>>>>>>>>>>>>>>>>>>>>>>>>>>>{\text{\Large \textcolor{white}{$\blacksquare$}}}
    |>>>>>>>>>>>>>>>>>>>>>>>>>{\text{\Large \textcolor{white}{$\blacksquare$}}}
    |>>>>>>>>>>>>>>>>>>>>>>{\text{\Large \textcolor{white}{$\blacksquare$}}}
& \ast & \ldots &  \ast & G_0 & \ast & \ldots & \ast \\ & \ast \ar[u] \ar[ur] & \ast \ar[u] \ar[ur] & \ldots & \ast \ar[u] \ar[ur] & \ast \ar[u]  \ar[ur] & \ast \ar[u] \ar[ur] & \ldots & \  \ast  , \ar[u] }$$

$$\xymatrix{A_1 \colon & \ast
    \ar@{<-}[drrrrrrr]%
    |<<<<<<<<<<<{\text{\Large \textcolor{white}{$\blacksquare$}}}%
    |<<<<<<<<<<<<<<<<<<<<<{\text{ \Large \textcolor{white}{$\blacksquare$}}}
    |<<<<<<<<<<<<<<<<<<<<<<<<<<<<<<<<<<<<<{\text{\Large \textcolor{white}{$\blacksquare$}}}
    |<<<<<<<<<<<<<<<<<<<<<<<<<<<<<<<<<<<<<<<<<<<<{\text{\Large \textcolor{white}{$\blacksquare$}}}
    |>>>>>>>>>>>>>>>>>>>>>>>>>>>>>>>>>>>>>>{\text{\Large \textcolor{white}{$\blacksquare$}}}
    |>>>>>>>>>>>>>>>>>>>>>>>>>>>>>>>>>{\text{\Large \textcolor{white}{$\blacksquare$}}}
    |>>>>>>>>>>>>>>>>>>>>>>>>>{\text{ \textcolor{white}{$\blacksquare$}}}
    |>>>>>>>>>>>>>>>>>>>>>>{\text{ \textcolor{white}{$\blacksquare$}}}
& \ast & \ldots &  \ast & G_1 & \ast & \ldots & \ast \\ & \ast \ar[u] \ar[ur] & \ast \ar[u] \ar[ur] & \ldots & \ast \ar[u] \ar[ur] & \ast \ar[u]  \ar[ur] & \ast \ar[u] \ar[ur] & \ldots & \ \ast, \ar[u]}$$

$$\xymatrix{A_2 \colon & \ast
   \ar@{<-}[drrrrrrr]%
    |<<<<<<<<<<<{\text{\Large \textcolor{white}{$\blacksquare$}}}%
    |<<<<<<<<<<<<<<<<<<<<<{\text{ \Large \textcolor{white}{$\blacksquare$}}}
    |<<<<<<<<<<<<<<<<<<<<<<<<<<<<<<<<<<<<<{\text{\Large \textcolor{white}{$\blacksquare$}}}
    |<<<<<<<<<<<<<<<<<<<<<<<<<<<<<<<<<<<<<<<<<<<<{\text{\Large \textcolor{white}{$\blacksquare$}}}
    |>>>>>>>>>>>>>>>>>>>>>>>>>>>>>>>>>>>>>>>>>>>>>>>>>>>>>{\text{\Large \textcolor{white}{$\blacksquare$}}}
    |>>>>>>>>>>>>>>>>>>>>>>>>>>>>>>>>>>>>>>>>>>>{\text{\Large \textcolor{white}{$\blacksquare$}}}
    |>>>>>>>>>>>>>>>>>>>>>>>>>>>>>>>>{\text{\Large \textcolor{white}{$\blacksquare$}}}
    |>>>>>>>>>>>>>>>>>>>>>>>>>>{\text{\Large \textcolor{white}{$\blacksquare$}}}
& \ast & \ldots &  \ast & G_1 & (I,0) \wedge G_0 & \ldots & \ast \\ & \ast \ar[u] \ar[ur] & \ast \ar[u] \ar[ur] & \ldots & \ast \ar[u] \ar[ur] & G_0 \ar[u]^-{\omega}  \ar[ur]^{\chi} & \ast \ar[u] \ar[ur] & \ldots& \ \ast, \ar[u]}$$
and

$$\xymatrix{A_3 \colon & \ast
\ar@{<-}[drrrrrrr]%
    |<<<<<<<<<<<{\text{\Large \textcolor{white}{$\blacksquare$}}}%
    |<<<<<<<<<<<<<<<<<<<<<{\text{ \Large \textcolor{white}{$\blacksquare$}}}
    |<<<<<<<<<<<<<<<<<<<<<<<<<<<<<<<<<<<<<{\text{\Large \textcolor{white}{$\blacksquare$}}}
    |<<<<<<<<<<<<<<<<<<<<<<<<<<<<<<<<<<<<<<<<<<<<<<<<<<<<<<<<<<<{\text{ \textcolor{white}{$\blacksquare$}}}
    |>>>>>>>>>>>>>>>>>>>>>>>>>>>>>>>>>>>>>>>>>>>>>>>>>>>>>{\text{\Large \textcolor{white}{$\blacksquare$}}}
    |>>>>>>>>>>>>>>>>>>>>>>>>>>>>>>>>>>>>>>>>>>>>>>>{\text{\Large \textcolor{white}{$\blacksquare$}}}
    |>>>>>>>>>>>>>>>>>>>>>>>>>>>>>>>>>>>>>>>>{\text{\Large \textcolor{white}{$\blacksquare$}}}
    |>>>>>>>>>>>>>>>>>>>>>>>>>>>>>>{\text{\Large \textcolor{white}{$\blacksquare$}}}
    |>>>>>>>>>>>>>>>>>>>>>>>>>{\text{\Large \textcolor{white}{$\blacksquare$}}}
& \ast & \ldots &  \ast & \ast & S^1 \wedge G_0 & \ldots & \ast \\ & \ast \ar[u] \ar[ur] & \ast \ar[u] \ar[ur] & \ldots & \ast \ar[u] \ar[ur] & \ast  \ar[u]  \ar[ur] & \ast \ar[u] \ar[ur] & \ldots& \ \ast, \ar[u]}$$
where $\chi \colon G_0 \longrightarrow (I,0) \wedge G_0$ is the canonical morphism and

$$(A_0)_{\zeta_i}=G_0,$$
$$(A_1)_{\zeta_i}=G_1,$$
$$(A_2)_{\zeta_i}=G_1, \;\;\; (A_2)_{\zeta_{i+1}}= (I,0) \wedge G_0,\;\;\; (A_2)_{\beta_i}= G_0,$$
$$(A_3)_{\zeta_{i+1}}=S^1 \wedge G_0.$$
Next, define the non-trivial entries of ${\alpha_0}$, ${\alpha_1}$, and ${\alpha_2}$ by
$$({\alpha_0})_{\zeta_i} = \omega, \;\;\; ({\alpha_1})_{\zeta_i} = 1, \;\;\;({\alpha_2})_{\zeta_{i+1}} = \pi \wedge 1,$$
where $\pi \colon I \longrightarrow S^1$ is the projection. It follows from the construction of $\Q$ that the triangle
$$\xymatrix{\Q(A_0) \ar[r]^{\Q(\alpha_0)} & \Q(A_1) \ar[r]^{\Q(\alpha_1)} & \Q(A_2) \ar[r]^-{\Q(\alpha_2)} & \Q(A_3) \cong \Q(A_0)[1]}$$
is isomorphic to (\ref{dim1tr}). (Note that the morphism $ \xi' \colon \Cone(\chi) \longrightarrow S^1 \wedge G_0$ induced from the arrow morphism
$$\xymatrix{ G_0 \ar[r]^-{\chi} \ar[d] & (I,0) \wedge G_0 \ar[d]^{\pi \wedge 1} \\ \ast \ar[r] & S^1 \wedge G_0}$$
on the cones and the morphism $ \xi \colon \Cone(\chi) \longrightarrow S^1 \wedge G_0$ coming from the cone construction of \ref{cone} differ by a sign in the homotopy category $\Ho(\M)$. This explains the appearance of the minus sign before $\partial$ in (\ref{dim1tr}).)

Applying the functor $\Hocolim$ to the sequence (\ref{dim1seq}) gives the elementary distinguished triangle
$$\xymatrix{G_0 \ar[r]^{\omega} & G_1 \ar[r] & \Cone(\omega) \ar[r] & S^1 \wedge G_0}$$
since $A_0$, $A_1$, $A_2$ and $A_3$ are cofibrant. Now, by construction of $\R$, the desired result follows. \end{prf}

\begin{lem} \label{shenaxva} Let $\M$, $P$ and $N$ be as in \ref{I}. Assume that the triangle
$$ \xymatrix{ F_1 \ar[r]^{\iota} & F_0 \ar[r] & M \ar[r] & F_1[1]}$$
is antidistinguished in $\D(\pi_*P)$, and $H_*F_0$, $H_*F_1$ are projective $\pi_*P$-modules, and $H_*(\iota)$ is a monomorphism. Then the triangle
$$ \xymatrix{ \R(F_1) \ar[r]^{\R(\iota)} & \R(F_0) \ar[r] & \R(M) \ar[r] & \R(F_1[1]) \cong \Sigma \R(F_1)}$$
is distinguished. \end{lem}

\begin{prf} By \ref{dashla}, we have splittings
$$F_0 \cong \oplus_{i \in \Z/N\Z} F_0^{(i)}\;\;\;\text{and} \;\;\;  F_1 \cong \oplus_{i \in \Z/N\Z} F_1^{(i)}$$
in $\D(\pi_*P)$, where $H_*F_0^{(i)}, H_*F_1^{(i)} \in \B[i]$. Since $H_*F_1$ is projective, by \ref{proeqciulebi} (i), it follows that for any $i,j \in \Z/N\Z$,
$$[F_1^{(i)}, F_0^{(j)}] \cong \Hom_{\pi_*P}(H_*F_1^{(i)}, H_*F_0^{(j)}).$$
In particular,
$$[F_1^{(i)}, F_0^{(j)}] =0,\;\;\;\; i \neq j,$$
whence,
$$[F_1, F_0] \cong \oplus_{i \in \Z/N\Z} [F_1^{(i)}, F_0^{(i)}].$$
This implies that the morphism $ \iota \colon F_1 \longrightarrow F_0$ splits as well. More precisely, there are morphisms $\iota^{(i)} \colon F_1^{(i)} \longrightarrow F_0^{(i)}$ such that the diagram
$$\xymatrix{ F_1 \ar[rr]^{\iota} \ar[d]_{\cong} & & F_0 \ar[d]^{\cong} \\ \oplus_i F_1^{(i)} \ar[rr]^{\oplus_i \iota^{(i)}} & & \oplus_i F_0^{(i)}}$$
commutes in $\D(\pi_*P)$. Using this we see that the antidistinguished triangle
$$ \xymatrix{ F_1 \ar[r]^{\iota} & F_0 \ar[r] & M \ar[r] & F_1[1]}$$
is isomorphic to a finite sum of triangles of the form (\ref{dim1tr}). As $\R$ is an additive functor, \ref{trsh} completes the proof. \end{prf}

\subsection{Proof of \ref{I} (ii)} \label{daboloeba1} First we show the essential surjectivity. Consider $X \in \Ho(\M)$ with $\proj  \pi_*X \leq 1$. Since
$$H_* \colon \D(\pi_*P) \longrightarrow \Grmod \pi_*P$$
is essentially surjective, there is $M \in \D(\pi_*P)$ such that $H_*M$ and $\pi_*X$ are isomorphic as graded $\pi_*P$-modules. By \ref{I} (i), $\pi_*\R \cong H_*$ and, therefore, $\pi_*(\R(M)) \cong \pi_*X$ in $\Grmod \pi_*P$. As $\proj  \pi_*X \leq 1$, it follows from \ref{uct} that
$$\R(M) \cong X.$$

Next let us verify fully faithfulness. Let $M, M' \in \D(\pi_*P)$ and suppose that $\proj  H_*M$ is at most one. Now we check that the morphism
$\xymatrix{\R \colon [M,M'] \ar[r] & [\R(M), \R(M')]}$
is an isomorphism. By \ref{proeqciulebi}, one can choose $F \in \D(\pi_*P)$ and
$$\xymatrix{F \ar[r]^{\sigma} & M}$$
such that $H_*F$ is a projective $\pi_*P$-module and $H_*(\sigma)$ is an epimorphism. Then, embed $\sigma$ into a distinguished triangle
$$\xymatrix{Y \ar[r] & F \ar[r]^{\sigma} & M \ar[r] & Y[1].}$$
Since $\sigma$ induces a surjection on $H_*$, we have a short exact sequence
$$\xymatrix{0 \ar[r] & H_*Y \ar[r] & H_*F \ar[r] & H_*M \ar[r] & 0.}$$
As $\proj  H_*M \leq 1$ and $H_*F$ is projective, $H_*Y$ is projective as well. Therefore, by \ref{shenaxva}, the triangle
$$\xymatrix{\R(Y) \ar[r] & \R(F) \ar[r]^-{\R(\sigma)} & \R(M) \ar[r] & \R(Y[1]) \cong \Sigma \R(Y)}$$
is antidistinguished. In particular, if we apply $[-, \R(M')]$ to this triangle, we get a long exact sequence. Moreover, $\R$ induces a morphism of long exact sequences
{\footnotesize
$$\xymatrix{[F[1],M'] \ar[r] \ar[d]^\R & [Y[1], M'] \ar[r] \ar[d]^\R & [M, M'] \ar[r] \ar[d]^\R & [F, M'] \ar[r] \ar[d]^\R & [Y, M'] \ar[d]^\R \\  [\R(F[1]), \R(M')] \ar[r] & [\R(Y[1]), \R(M')] \ar[r] & [\R(M), \R(M')]  \ar[r] & [\R(F), \R(M')] \ar[r] & [\R(Y), \R(M')].}$$}

\noindent{}By \ref{I} (i), Lemma \ref{proeqciulebi} (i) and the five lemma,
$$\xymatrix{\R \colon [M,M'] \ar[r] & [\R(M), \R(M')] }$$
is an isomorphism. $\Box$

\begin{coro}\label{dim1case} Let $\M$ be a stable model category and $P$ a compact generator of $\Ho(\M)$. Suppose that the graded ring $\pi_*P$ is $N$-sparse and $\gl \pi_*P = 1.$ Then the functor
$$\R \colon \D(\pi_*P) \longrightarrow \Ho(\M)$$
is an equivalence of categories. \end{coro}

\begin{remk} Lemma \ref{shenaxva} and the construction of Adams spectral sequence \ref{ASS} show that $\R$ preserves (up to sign) universal coefficient sequences. More precisely, for any $M \in \D(\pi_*P)$ with $\proj  H_*M \leq 1$ and $M' \in \D(\pi_*P)$, the diagram
{\scriptsize
$$ \xymatrix{  0 \ar[r] & \Ext^{1}_{\pi_*P}( H_*M[1], H_*M') \ar[r] \ar[d]^{\cong} & [M,M'] \ar[r]^-{H_*} \ar[d]^{\R} & \Hom_{\pi*(P)} (H_*M, H_*M') \ar[r] \ar[d]^{\cong} & 0 \\  0 \ar[r] & \Ext^{1}_{\pi_*P}( \pi_*(\R(M))[1], \pi_*(\R(M'))) \ar[r] & [\R(M),\R(M')] \ar[r]^-{\pi_*} & \Hom_{\pi_*P} (\pi_*(\R(M)), \pi_*(\R(M'))) \ar[r] & 0}$$}
commutes up to sign, where the left and right isomorphisms come from the functor isomorphism $$\pi_* \circ \R \cong H_*.$$ \end{remk}

Note that the equivalence of categories
$$\D(\pi_*P) \sim \Ho(\M)$$
in \ref{dim1case} is a special case of \cite[4.3.2]{G99}. Moreover, since the functor $\R$ preserves universal coefficient sequences, it is in fact naturally isomorphic to the functor constructed in \cite[4.3]{G99}. One should note that although \cite[4.3.2]{G99} is stated for homotopy categories coming from stable model categories, the proof in \cite[4.3]{G99} can be equally applied to general triangulated categories as well. In other words, one has the following generalization of \ref{dim1case}.

\begin{prop} \label{nontreq} Suppose $\T$ is a triangulated category with infinite coproducts, $P \in \T$ a compact generator, and let $\pi_* = \Hom_{\T}(P,-)_*$ Suppose that the graded ring $\pi_*P$ is $N$-sparse and $\gl\pi_*P = 1$. Then $\T$ is equivalent to $\D(\pi_*P)$. \end{prop}

\begin{remk} According to \cite[4.3.1]{G99} the equivalence $\D(\pi_*P) \sim \Ho(\M)$ is even a triangulated equivalence. However the present author was unable to follow the argument about Toda brackets at the end of \cite[4.3]{G99}. We do not know whether the equivalences \ref{nontreq} and \ref{dim1case} are triangulated or not. \end{remk}

\subsection{Examples}

\begin{exmp} [Complex $K$-theory] It is well known to specialists that the derived category $\D(KU)$ of the complex $K$-theory spectrum $KU$ is equivalent to the derived category $\D(\pi_*KU)$ of the homotopy ring
$$\pi_*KU \cong \Z[u, u^{-1}],\;\;\;\;\;\;\;\;\;\;\;\; |u|=2.$$
Since $\pi_*KU$ is $2$-sparse and $\gl \pi_*KU=1$, one can think of this equivalence as a corollary of \ref{dim1case}. Note also that $\Omega^{\infty}KU$ is not a product of Eilenberg-MacLane spaces. Consequently, in view of \ref{noneq2}, there does not exist a zig-zag of Quillen equivalences between the model categories $\Mod KU$ and $\Mod \Z[u, u^{-1}]$. In particular,
$$\R \colon \D(\Z[u, u^{-1}]) \longrightarrow \D(KU)$$
can not be derived from a zig-zag of Quillen equivalences. Note that we do not know whether the equivalence $\D(\Z[u, u^{-1}]) \sim \D(KU)$ is triangulated or not.
\end{exmp}

\begin{exmp} [Connective Morava $K$-theories] The connective Morava $K$-theory spectrum $k(n)$, $n\geq 1$, is obtained from the truncated Brown-Peterson spectrum $BP\langle n\rangle$ by killing the regular sequence $(p, v_1,\ldots,v_{n-1})$. In particular, we have an isomorphism of graded rings
$$\pi_*k(n) \cong \mathbb{F}_{p}[v_n],\;\;\;\;\;\;\;\;\;\;\;\; |v_n|= 2(p^n-1).$$
As $\gl\pi_*k(n)=1$ and $\pi_*k(n)$ is $2(p^n-1)$-sparse, the functor
$$\R \colon \D(\mathbb{F}_{p}[v_n]) \longrightarrow \D(k(n))$$
is an equivalence of categories by \ref{dim1case}. Observe that $\Omega^{\infty}k(n)$ is not a product of Eilenberg-MacLane spaces. Therefore, by \ref{noneq2}, the model categories $\Mod k(n)$ and $\Mod \mathbb{F}_{p}[v_n]$ are not connected by a zig-zag of Quillen equivalences. In particular, $\R$ does not come from a zig-zag of Quillen equivalences.  \end{exmp}

\section{The two dimensional case}\label{intgamoy}

In this section we prove \ref{I} (iii) and \ref{equival}.

\subsection{Proof of \ref{I} (iii)}

Let $X \in \Ho(\M)$ and suppose $\proj  \pi_*X \leq 2$. By \ref{proeqciulebi} and \ref{dashla}, there are cofibrant objects $F^{(i)}$, $i \in \Z/ N\Z$, and a morphism
$$\xymatrix{\vee_i F^{(i)} \ar[r]^-{\sigma} & X}$$
such that $\pi_*F^{(i)}$ are projective, $\pi_*F^{(i)} \in \B[i]$ and
$$\xymatrix{\oplus_i \pi_*F^{(i)} \ar[r]^-{\pi_*\sigma} & \pi_*X}$$
is an epimorphism. Embed $\sigma$ into a distinguished triangle
\begin{align}\label{samkutxedi1} \xymatrix{Y \ar[r] & \vee_i F^{(i)} \ar[r]^-{\sigma} & X \ar[r] & \Sigma Y.}\end{align}
Since $\sigma$ induces a surjection on $\pi_*$, we have a short exact sequence
$$\xymatrix{0 \ar[r] & \pi_*Y \ar[r] & \oplus_i \pi_*F^{(i)} \ar[r]^-{\pi_*\sigma} & \pi_*X \ar[r] & 0.}$$
As $\proj  \pi_*X \leq 2$ and $\oplus_i \pi_*F^{(i)}$ is projective, one has
$$\proj  \pi_*Y \leq 1.$$
Combining this with \ref{dashla}, we get a splitting
$$Y \cong \vee_i Y^{(i)}$$
in $\Ho(\M)$, where $Y^{(i)} \in \M_{\cof}$, and $\pi_*Y^{(i)} \in \B[i]$, $i \in \Z/ N\Z$. If we replace $Y$ by $\vee_i Y^{(i)}$ in (\ref{samkutxedi1}), we get a distinguished triangle
$$\xymatrix{\vee_i Y^{(i)} \ar[r] & \vee_i F^{(i)} \ar[r]^{\sigma} & X \ar[r] & \Sigma (\vee_i Y^{(i)}).}$$
By Adams spectral sequence \ref{ASS},
$$[Y^{(i)}, F^{(j)}] = 0,\;\;\;\;\;\; j \neq i, i+1 \; (\text{mod} \; N),$$
since $\proj  \pi_*Y^{(i)} \leq 1.$ Therefore the only possible nontrivial entries of $$ \vee_i Y^{(i)} \longrightarrow \vee_i F^{(i)}$$ are
$$\xymatrix{Y^{(i)} \ar[r]^{\alpha_i} & F^{(i)}}\;\;\;\; \text{and} \;\;\;\; \xymatrix{Y^{(i-1)} \ar[r]^{\beta_i} & F^{(i)}.}$$
Without loss of generality one may assume that $F^{(i)}$ is bifibrant for any $i \in \Z/N\Z$. Hence there exist morphisms
$$ l_i \colon Y^{(i)} \longrightarrow F^{(i)}\;\;\;\; \text{and} \;\;\;\; k_i \colon Y^{(i-1)} \longrightarrow F^{(i)}$$
in $\M$ whose homotopy classes are $\alpha_i$ and $- \beta_i$, respectively. Consequently, we get a diagram
$$\xymatrix{W \colon & F^{(0)}
    \ar@{<-}[drrrrrrrr]%
    |<<<<<<<<<<<<<<<<<<{\text{\Large \textcolor{white}{$\blacksquare$}}}%
    |<<<<<<<<<<<<<<<<<<<<<<<{\text{\Large \textcolor{white}{$\blacksquare$}}}
    |<<<<<<<<<<<<<<<<<<<<<<<<<<<<<<<<<<<<<<<{\text{\Large \textcolor{white}{$\blacksquare$}}}
    |>>>>>>>>>>>>>>>>>>>>>>>>>>>>>>>>>>>>>>>>>{\text{\Large \textcolor{white}{$\blacksquare$}}}
    |>>>>>>>>>>>>>>>>>>>>>>>>>>{\text{\Large \textcolor{white}{$\blacksquare$}}}
    |>>>>>>>>>>>>>>>>>>{\text{\Large \textcolor{white}{$\blacksquare$}}}
    & & F^{(1)} & & \ldots & & F^{(N-2)} & & F^{(N-1)} \\  & Y^{(0)} \ar[u]^{l_0} \ar[urr]^{k_1} & & Y^{(1)} \ar[u]_{l_1} \ar[urr] & & \ldots \ar[urr] & & Y^{(N-2)} \ar[u]^{l_{N-2}} \ar[urr]^{k_{N-1}} & &  Y^{(N-1)}, \ar[u]_{l_{N-1}}}$$
in $\M$, i.e., an object of $\M^{\C_N}$. We may also assume that $W \in (\M^{\C_N})_{\cof}$ (otherwise we could replace it cofibrantly). The homotopy colimit of $W$ is isomorphic to the homotopy coequalizer of
$$\xymatrix{\vee_{i}Y^{(i)} \ar@<0.5ex>[r]^l \ar@<-0.5ex>[r]_k & \vee_{i}F^{(i)}},$$
where $k= \vee_i k_i$ and $l=\vee_i l_i$. Therefore, by \ref{htpycoeqtr}, one has
$$\Hocolim W \cong X$$
in $\Ho(\M)$. On the other hand, clearly, $W \in \K$. Hence the construction of $\R$ implies
$$\R(\Q(W)) \cong X.\;\;\;\;\;\; \Box$$

\subsection{Technical lemmas}

In this subsection we prove two lemmas which will be used in the proof of \ref{equival}.

\begin{lem}\label{techlemma} Let $\M$, $P$, $\pi_*$ and $N$ be as in \ref{equival}, and $(C_*, d)$ be a DG $\pi_*P$-module and $G_*$ a DG $\pi_*P$-module with zero differential (i.e., a graded $\pi_*P$-module), and suppose that
$$\xymatrix{f \colon G_* \longrightarrow C_*}$$
is a DG morphism which induces a monomorphism on homology. Then the functor
$$\R \colon \D(\pi_*P) \longrightarrow \Ho(\M)$$
sends the antidistinguished triangle (see \ref{algcone})
\begin{align}\label{dreieck}\xymatrix{G_* \ar[r]^f & C_* \ar[r]^-{\iota} & C(f) \ar[r]^{-\partial} & G_*[1]}\end{align}
to a distinguished triangle.
\end{lem}

\begin{prf} The proof goes in the same way as that of \ref{trsh}. First we construct a sequence (see \ref{commshift})
$$\xymatrix{X \ar[r] & Y \ar[r] & Z \ar[r] & X^\# }$$
in $\L$ ($\L$=$\K$ under the hypotheses of \ref{equival}) which is mapped to a triangle isomorphic to (\ref{dreieck}) after applying $\Q$.

Consider the splitting
$$G_* \cong G^{(0)} \oplus \ldots \oplus G^{(N-1)},$$
in $\Grmod \pi_*P$, where $G^{(i)} \in \B[i]$, $i \in \Z/N \Z$. In view of \ref{II} (i), there are cofibrant $X_{\zeta_i}$-s with
$$\pi_*X_{\zeta_i} \cong G^{(i)}.$$
Clearly, the diagram
$$\xymatrix{X \colon & X_{\zeta_0}
    \ar@{<-}[drrrr]%
    |<<<<<<<<<{\text{\Large \textcolor{white}{$\blacksquare$}}}%
    |<<<<<<<<<<<<<{\text{\Large \textcolor{white}{$\blacksquare$}}}
    |<<<<<<<<<<<<<<<<<<<<<<{\text{\Large \textcolor{white}{$\blacksquare$}}}
    |>>>>>>>>>>>>>>>>>>>>>>>>>>>>>{\text{\Large \textcolor{white}{$\blacksquare$}}}
    |>>>>>>>>>>>>>>>>>>>{\text{\Large \textcolor{white}{$\blacksquare$}}}
    |>>>>>>>>>>>>>>{\text{\Large \textcolor{white}{$\blacksquare$}}}
& X_{\zeta_1} & \ldots & X_{\zeta_{N-2}} & X_{\zeta_{N-1}} \\ & \ast \ar[u] \ar[ur] & \ast \ar[u] \ar[ur] & \ldots \ar[ur] & \ast \ar[u] \ar[ur] & \ast \ar[u]}$$
is an object of $\L$ and $\Q(X) \cong (G_*, 0)$. Next, by \ref{Qequicase}, $\Q(\L) = \Mod \pi_*P$ since $\gl \pi_*P$ is less than $N-1$. In particular, one can choose a bifibrant
$$\xymatrix{Y \colon & Y_{\zeta_0}
    \ar@{<-}[drrrrrrrr]%
    |<<<<<<<<<<<<<<<<<<{\text{\Large \textcolor{white}{$\blacksquare$}}}%
    |<<<<<<<<<<<<<<<<<<<<<<<{\text{\Large \textcolor{white}{$\blacksquare$}}}
    |<<<<<<<<<<<<<<<<<<<<<<<<<<<<<<<<<<<<<<{\text{\Large \textcolor{white}{$\blacksquare$}}}
    |>>>>>>>>>>>>>>>>>>>>>>>>>>>>>>>>>>>>>>>>{\text{\Large \textcolor{white}{$\blacksquare$}}}
    |>>>>>>>>>>>>>>>>>>>>>>>>{\text{\Large \textcolor{white}{$\blacksquare$}}}
    |>>>>>>>>>>>>>>>>>>{\text{\Large \textcolor{white}{$\blacksquare$}}}
& & Y_{\zeta_1} & & \ldots & & Y_{\zeta_{N-2}} & & Y_{\zeta_{N-1}} \\  & Y_{\beta_0} \ar[u]^{l_0} \ar[urr]^{k_1} & & Y_{\beta_1} \ar[u]_{l_1} \ar[urr] & & \ldots \ar[urr] & & Y_{\beta_{N-2}} \ar[u]^{l_{N-2}} \ar[urr]^{k_{N-1}} & & Y_{\beta_{N-1}}, \ar[u]_{l_{N-1}}}$$
in $\L$ such that $$\Q(Y) \cong (C_*,d).$$ It then follows from \ref{Qequicase} that there is a morphism $\psi \colon X \longrightarrow Y$ in $\L$ so that the diagram
$$\xymatrix{ \Q(X) \ar[r]^{\Q(\psi)} \ar[d]^{\cong} & \Q(Y) \ar[d]^{\cong} \\ G_* \ar[r]^f & C_*}$$
commutes in $\Mod \pi_*P$. Further, as $X$ is cofibrant and $Y$ is fibrant, the morphism $\psi$ comes from some morphism in $\M^{\C_N}$, denoted by $\psi$ as well. Moreover, using the factorization axiom for the model category $\M^{\C_N}$, we may assume that $\psi \colon X \longrightarrow Y$ is cofibration $\M^{\C_N}$, i.e.,
$$\xymatrix{X_{\zeta_i} \vee Y_{\beta_{i-1}} \vee Y_{{\beta}_i} \ar[rr]^-{(\psi_{\zeta_i}, k_{i-1}, l_i)} & &   Y_{{\zeta}_i}}$$
is a cofibration for any $i \in \Z / N \Z$ (see \ref{Reedy}). Next, we construct an object in $\L$ whose image under $\Q$ is isomorphic to the algebraic mapping cone $C(f)$. Define $Z \in \M^{\C_N}$ by
$$Z_{\zeta_i}= CX_{\zeta_{i-1}} \vee Y_{\zeta_i},\;\;\;\;\;\; Z_{\beta_i}=X_{\zeta_i} \vee Y_{\beta_i},$$
$$\xymatrix{\lambda_i \colon Z_{\beta_i} =X_{\zeta_i} \vee Y_{\beta_i} \ar[rr]^-{(\psi_{\zeta_i}, l_i)} & & Y_{\zeta_i}  \ar[rr]^-{\incl} & & CX_{\zeta_{i-1}} \vee Y_{\zeta_i}=Z_{\zeta_i},}$$
$$\xymatrix{\kappa_i \colon Z_{\beta_{i-1}}= X_{\zeta_{i-1}} \vee Y_{\beta_{i-1}} \ar[rr]^-{\varepsilon_{i-1} \vee k_i} & & CX_{\zeta_{i-1}} \vee Y_{\zeta_i}=Z_{\zeta_i},}$$
where $$C(-)= (I,0) \wedge -,$$ and $\varepsilon_i \colon X_{\zeta_i} \longrightarrow CX_{\zeta_i}$, $i \in \Z /N\Z$, is the canonical map. Thus, the diagram $Z$ looks as follows:
{\footnotesize
$$\hspace{-1cm}\xymatrix{CX_{\zeta_{-1}} \vee Y_{\zeta_0}
\ar@{<-}[drrrrrrrr]%
    |<<<<<<<<<<<<<<<<<<<<<<<{\text{\Large \textcolor{white}{$\blacksquare$}}}%
    |<<<<<<<<<<<<<<<<<<<<<<<<<<<<<<<<<{\text{\Large \textcolor{white}{$\blacksquare$}}}
    |<<<<<<<<<<<<<<<<<<<<<<<<<<<<<<<<<<<<<<<<<<<<<<<<<<<<<{\text{\Large \textcolor{white}{$\blacksquare$}}}
    |>>>>>>>>>>>>>>>>>>>>>>>>>>>>>>>>>>>>>>>>>>>>>>>>>>>>>>>>>>>{\text{\Large \textcolor{white}{$\blacksquare$}}}
    |>>>>>>>>>>>>>>>>>>>>>>>>>>>>>>>>>>>>>>{\text{\Large \textcolor{white}{$\blacksquare$}}}
    |>>>>>>>>>>>>>>>>>>>>>>>>>>{\text{\Large \textcolor{white}{$\blacksquare$}}}
& &CX_{\zeta_0} \vee Y_{\zeta_1} & & \ldots & & CX_{\zeta_{N-3}} \vee Y_{\zeta_{N-2}}& & CX_{\zeta_{N-2}} \vee Y_{\zeta_{N-1}} \\   X_{\zeta_0} \vee Y_{\beta_0} \ar[u]^{\lambda_0} \ar[urr]^{\kappa_1} & & X_{\zeta_1} \vee Y_{\beta_1}  \ar[u]_{\lambda_1} \ar[urr] & & \ldots \ar[urr] & & X_{\zeta_{N-2}} \vee Y_{\beta_{N-2}}  \ar[u]^{\lambda_{N-2}} \ar[urr]^{\kappa_{N-1}} & & X_{\zeta_{N-1}} \vee Y_{\beta_{N-1}}. \ar[u]_{\lambda_{N-1}}}$$}

\noindent{}As $X$ and $Y$ are cofibrant and $\psi \colon X \longrightarrow Y$ is a cofibration, it follows that $Z \in (\M^{\C_N})_{\cof}$. Further,
$$\xymatrix{\pi_*\lambda_i \colon \pi_*Z_{\beta_i} \ar[r] & \pi_*Z_{\zeta_i}}$$
is a monomorphism since $f$ induces a monomorphism on homology. Besides, obviously, $\pi_*Z_{\beta_i}$ and $\pi_*Z_{\zeta_i}$ are objects of $\B[i]$. Consequently, one gets $Z \in \L$. The construction of $\Q$ immediately implies that
$$\Q(Z) \cong C(f)$$
in $\Mod \pi_*P$.

The diagram $Z$ comes with two obvious morphisms $Y \longrightarrow Z$ and $Z \longrightarrow X^\#$ in $\M^{\C_N}$. Hence we get a sequence
\begin{align}\label{shualeduri}\xymatrix{X \ar[r]^\psi & Y \ar[r] & Z \ar[r] & X^\#}\end{align}
in $\L$. It is immediate from the construction of $\Q$ that the triangle
$$\xymatrix{\Q(X) \ar[r]^{\Q(\psi)} & \Q(Y) \ar[r] & \Q(Z) \ar[r] & \Q(X^\#) \cong \Q(X)[1]}$$
is isomorphic to (\ref{dreieck}) (the sign appears by the same reason as in \ref{trsh}). On the other hand, applying the functor $\Hocolim$ to (\ref{shualeduri}) ($X$, $Y$, $Z$ are cofibrant) one obtains the elementary distinguished triangle
$$\xymatrix{\colim X \ar[r]^{\colim \psi} & \colim Y \ar[r] & \Cone( \colim \psi) \ar[r]^\delta & S^1 \wedge \colim X.}$$
Now the desired result follows from the construction of $\R$. \end{prf}

\begin{lem}\label{dreieck1} Let $\M$, $P$, $N$ be as in \ref{equival}. Suppose that the triangle
\begin{align}\label{dreieck2}\xymatrix{M' \ar[r]^{\iota} & F \ar[r]^f & M \ar[r]^-{-d} & M'[1]}\end{align}
is antidistinguished in $\D(\pi_*P)$, and assume that $\proj  H_*M' \leq 1$ and $H_*(\iota)$ is a monomorphism. Then the triangle
$$\xymatrix{\R(M') \ar[r]^{\R(\iota)} & \R(F) \ar[r]^{\R(f)} & \R(M) \ar[r]^-{-\R(d)} & \R(M'[1]) \cong \Sigma \R(M')}$$
is distinguished in $\Ho(\M)$. \end{lem}

\begin{prf} By \ref{uct}, since $\proj  H_*M' \leq 1$, one has
$$M' \cong H_*M'$$
in $\D(\pi_*P)$, where $H_*M'$ is regarded as a differential graded $\pi_*P$-module with zero differential. Therefore, the triangle (\ref{dreieck2}) is isomorphic to an antidistinguished triangle
$$\xymatrix{H_*M' \ar[r] & F \ar[r]^f & M \ar[r] & H_*M'[1].}$$
Without loss of generality we may assume that $M$ is cofibrant in $\Mod \pi_*P$. Then the last morphism in the latter triangle comes from some morphism $\varphi \colon M \longrightarrow H_*M'[1]$ of differential graded modules and, by the triangulated category axiom \textbf{(T3)}, one has a triangle isomorphism
$$\xymatrix{H_*M' \ar[r] \ar@{=}[d] & F \ar[r] \ar@{-->}[d]^{\cong} & M \ar[r]^-{\varphi} \ar@{=}[d] & H_*M'[1] \ar@{=}[d] \\
H_*M' \ar[r]^-{\incl} & C(\varphi)[-1] \ar[r]^-{-\partial[-1]}  & M \ar[r]^-{\varphi} & H_*M'[1] }$$
in $\D(\pi_*P)$ (see \ref{algcone}). On the other hand, the same axiom shows that the lower triangle is ismorphic to
$$\xymatrix{H_*M' \ar[r]^-{\incl} & C(\varphi)[-1] \ar[r] & C(\incl) \ar[r]^{-\partial} & H_*M'[1].}$$
Thus the triangle (\ref{dreieck2}) is isomorphic to an antidistinguished triangle of the form (\ref{dreieck}) from \ref{techlemma}. This completes the proof.
\end{prf}

\subsection{Proof of \ref{equival}} \label{daboloeba2} It remains to check that $\R$ is fully faithful when restricted to the full subcategory of at most two dimensional objects. We now verify a more general fact. Namely, we show that for $M \in \D(\pi_*P)$ with $\proj  H_*M \leq 2$ and any $M' \in \D(\pi_*P)$, the morphism
$$\xymatrix{\R \colon [M,M'] \ar[r] & [\R(M), \R(M')] }$$
is an isomorphism. Indeed, by \ref{proeqciulebi}, we may choose $F \in \D(\pi_*P)$ and $\sigma \colon F \longrightarrow M$ so that $H_*F$ is projective and $H_*(\sigma)$ an epimorphism. Embed $\sigma$ into a distinguished triangle
$$\xymatrix{Y \ar[r] & F \ar[r]^{\sigma} & M \ar[r] & Y[1].}$$
As $\sigma$ induces a surjection on $H_*$, one has a short exact sequence
$$\xymatrix{0 \ar[r] & H_*Y \ar[r] & H_*F \ar[r]^{H_*(\sigma)} & H_*M \ar[r] & 0.}$$
It then follows that
$$\proj  H_*Y \leq 1$$
since $\proj  H_*M \leq 2$ and $H_*F$ is projective. Therefore, by \ref{dreieck1}, the triangle
$$\xymatrix{\R(Y) \ar[r] & \R(F) \ar[r]^{\R(\sigma)} & \R(M) \ar[r] & \R(Y[1]) \cong \Sigma \R(Y) }$$
is antidistinguished. In particular, if we apply $[-, \R(M')]$ to this triangle, we get a long exact sequence. Moreover, $\R$ induces a morphism of long exact sequences
{\footnotesize
$$\xymatrix{[F[1],M'] \ar[r] \ar[d]^\R & [Y[1], M'] \ar[r] \ar[d]^\R & [M, M'] \ar[r] \ar[d]^\R & [F, M'] \ar[r] \ar[d]^\R & [Y, M'] \ar[d]^\R \\  [\R(F[1]), \R(M')] \ar[r] & [\R(Y[1]), \R(M')] \ar[r] & [\R(M), \R(M')]  \ar[r] & [\R(F), \R(M')] \ar[r] & [\R(Y), \R(M')].}$$}

\noindent{}By \ref{proeqciulebi} (i), the first and fourth morphisms in the diagram are isomorphisms. Besides, the proof in Subsection \ref{daboloeba1} shows that so are the second and last one. Consequently, by the five lemma,
$$\xymatrix{\R \colon [M,M'] \ar[r] & [\R(M), \R(M')] }$$
is an isomorphism. $\Box$

\begin{remk} \label{shecdoma} Let $\M$, $P$ and $N$ be as in \ref{I}. Franke in \cite[2.2, Proposition 2]{F96} claims that for any distinguished triangle
$$\xymatrix{M \ar[r]^f & M' \ar[r]^g & M'' \ar[r]^{\partial} & M[1]}$$
in $\D(\pi_*P)$ with $H_*(f)$ a monomorphism, the triangle
$$\xymatrix{\R(M) \ar[r]^{\R(f)} & \R(M') \ar[r]^{\R(g)} & \R(M'') \ar[r]^-{\R({\partial})} & \R(M[1]) \cong \Sigma \R(M)}$$
is (anti)distinguished in $\Ho(\M)$. If this was so, then one could prove (similarly to the proofs of \ref{I} (ii) and \ref{equival}) that under the hypotheses of \ref{I}
$$\xymatrix{ \R \colon \D(\pi_*P) \ar[r] & \Ho(\M)}$$
is fully faithful. Moreover, using the claim, one can also check that $\R$ is essentially surjective and thus an equivalence of categories. However, as far as the present author can see, the proof of \cite[2.2, Proposition 2]{F96} seems to contain a gap. The point is that the arguments given in \cite[2.2]{F96}, show that the sequence
$$\xymatrix{\R(M) \ar[r]^{\R(f)} & \R(M') \ar[r]^{\R(g)} & \R(M'')}$$
is a part of some distinguished triangle, but do not say anything about why exactly the morphism
$$\xymatrix{\R(M'') \ar[r]^-{\R({\partial})} & \R(M[1]) \cong \Sigma \R(M)}$$
makes this sequence into an (anti)distinguished triangle.

We were unable to fill this gap in the general setting. Although we managed to do this in low dimensional cases and thus obtained the statements of \ref{I} and \ref{equival}.

Note that in particular, we do not know whether the functor $$\xymatrix{ \R \colon \D(\pi_*P) \ar[r] & \Ho(\M)}$$ is triangulated in general or not. \end{remk}

\subsection{Examples}

\begin{exmp}[Real connective $K$-theory localized at an odd prime]\label{ko_(p)} Consider $ko_{(p)}$, the real connective $K$-theory spectrum localized at an odd prime $p$. The ring isomorphism
$$\pi_*ko_{(p)} \cong \Z_{(p)}[\omega],\;\;\;\;\;\;\;\;\;\;\;\; |\omega|=4,$$
implies that $\gl \pi_*ko_{(p)}=2$ and $\pi_*ko_{(p)}$ is $4$-sparse. Thus \ref{equival2} can be applied to $ko_{(p)}$ and therefore we conclude that the functor
$$\xymatrix{ \R \colon \D(\Z_{(p)}[\omega]) \ar[r] & \D(ko_{(p)}) }$$
is an equivalence of categories. Note that $\Omega^{\infty}ko_{(p)}$ is not a product of Eilenberg-MacLane spaces. Hence, in view of \ref{noneq2}, there does not exist a zig-zag of Quillen equivalences between the model categories $\Mod ko_{(p)}$ and $\Mod \Z_{(p)}[\omega]$. In particular, the functor $$\xymatrix{ \R \colon \D(\Z_{(p)}[\omega]) \ar[r] & \D(ko_{(p)}) }$$ can not be derived from a zig-zag of Quillen equivalences.  \end{exmp}

\begin{exmp} The truncated Brown-Peterson spectrum $BP\langle 1 \rangle$ for an odd prime $p$ satisfies the assumptions of \ref{equival2}. Indeed, there is a ring isomorphism
$$\pi_*BP\langle 1 \rangle \cong \Z_{(p)}[v_1], \;\;\;\;\;\;\;\;\;\; |v_1|=2(p-1),$$
and $\gl \Z_{(p)}[v_1] =2$ and $2(p-1) \geq 4$. Thus, the functor
$$\R \colon \D(\Z_{(p)}[v_1]) \longrightarrow \D(BP \langle 1 \rangle)$$
is an equivalence of categories. Next, as mentioned above, $\Omega^{\infty}BP\langle 1 \rangle$ is not a product of Eilenberg-MacLane spaces. Therefore, by
\ref{noneq2}, the model categories $\Mod BP \langle 1 \rangle$ and $\Mod \Z_{(p)}[v_1]$  can not be connected by a zig-zag of Quillen equivalences. In particular, $\R \colon \D(\Z_{(p)}[v_1]) \longrightarrow \D(BP \langle 1 \rangle)$ is not derived from a zig-zag of Quillen equivalences. \end{exmp}

\begin{exmp} Another ring spectrum to which \ref{equival2} applies is the Johnson-Wilson spectrum $E(2)$ for an odd prime $p$. This follows from the ring isomorphism
$$\pi_*E(2) \cong \Z_{(p)}[v_1, v_2, v_2^{-1}], \;\;\;\;\;\;\;\;\;\; |v_1|= 2(p-1), \;\; |v_2| = 2(p^2-1),$$
since $\gl \Z_{(p)}[v_1, v_2, v_2^{-1}] =2$ and $2(p-1) \geq 4$. Further, as in the two previous examples, $\Omega^{\infty} E(2)$ is not a product of Eilenberg-MacLane spaces and therefore, by \ref{noneq2}, there is no zig-zag of Quillen equivalences between the model categories $\Mod E(2)$ and $\Mod \Z_{(p)}[v_1, v_2, v_2^{-1}]$. In particular, the equivalence $\R \colon \D(\Z_{(p)}[v_1, v_2, v_2^{-1}]) \longrightarrow \D(E(2))$ does not come from a zig-zag of Quillen equivalences. \end{exmp}

Finally, we would like to conclude the paper with the following

\begin{remk}\label{algebra} As noted above, \ref{dgcoro} and \ref{dgequival1} are corollaries of \ref{Rcoro} and \ref{equival1}, respectively. However, one can also give independent and purely algebraic proofs of \ref{dgcoro} and \ref{dgequival1} as well. The main idea is the same as in the proofs of \ref{I} and \ref{equival}, but in this special case one can use various algebraic constructions which are easier (compared to their topological counterparts) to deal with. Indeed, instead of the cone construction \ref{cone} we can apply the algebraic mapping cone construction of \ref{algcone}. Further, the mapping spaces in (\ref{plbk}) can be replaced by the internal Hom-complexes. Next, since any cofibration in $\Mod A$ is a split monomorphism of underlying graded modules, the right vertical arrow in (\ref{plbk}) becomes an epimorphism and thus (\ref{plbk}) leads to a short exact sequence of chain complexes. Clearly, the resulting long exact homology sequence plays the role of the Mayer-Vietoris sequence. Then, the verification of the identity $q\partial=b$ in (\ref{dididiag}) can be done by explicit algebraic calculations with homology classes. Finally, the long exact sequence of the homotopy coequalizer (\ref{hocoeq})  becomes the long exact homology sequence of
$$\xymatrix{0 \ar[r] & \oplus_{i}X_{\beta_{i}} \ar[r]^{l-k} & \oplus_{i}X_{\zeta_i} \ar[r] & \colim X \ar[r] & 0.} $$

The other points of the proofs of \ref{I} and \ref{equival} are more or less formal and can be directly used in the algebraic setting as well. \end{remk}

\begin{tabular}{l}
{\sc Mathematisches Institut}\\
{\sc Universit\"at Bonn}\\
{\sc Endenicher Allee 60}\\
{\sc 53115 Bonn, Germany}\\
\\[-10pt]
\emph{E-mail address}: irpatchk@math.uni-bonn.de
\end{tabular}

\end{document}